\documentclass[preprint,10pt, sort&compress]{elsarticle}

\usepackage[margin=2.5cm]{geometry}
\usepackage{setspace}
\usepackage{lmodern}
\usepackage{amsmath,amssymb}
\usepackage{amsfonts,amsthm}
\usepackage{lineno}
\usepackage{graphicx}
\usepackage{bm, bbm}
\usepackage{algorithm}
\usepackage{algpseudocode}

\usepackage{physics}
\usepackage{color}
\usepackage{pxfonts}
\usepackage[footnotesize,bf]{caption}
\usepackage{subcaption}
\usepackage{mathtools}
\usepackage{hhline}
\usepackage[T1]{fontenc}

\usepackage{caption}
\usepackage{graphicx}
\usepackage{float}
\usepackage{subcaption}
\usepackage{booktabs}

\usepackage{mathrsfs}
\usepackage{placeins}

\newtheorem{proposition}{Proposition}

\newtheorem{remark}{Remark}[section]

\numberwithin{equation}{section}
\numberwithin{figure}{section}
\numberwithin{table}{section}

\makeatletter
\def\ps@pprintTitle{%
	\let\@oddhead\@empty
	\let\@evenhead\@empty
	\let\@oddfoot\@empty
	\let\@evenfoot\@oddfoot
}
\makeatother

\graphicspath{{figures/}}
\usepackage{multirow}
\usepackage{ulem}
\usepackage{tikz}
\usepackage{color}
\usepackage[]{xcolor}

\begin{document}
\begin{frontmatter}
\title{Structure-Preserving Neural Time-Stepping for Gradient-Flow PDEs via Onsager's Variational Principle}

\author[a]{Zhonghua Qiao}
\author[b]{Xianmin Xu}
\author[a]{Rukang You}
\author[b]{Tao Zhou}

\address[a]{Department of Applied Mathematics, the Hong Kong Polytechnic University, Hung Hom, Hong Kong}
\address[b]{SKLMS, Academy of Mathematics and Systems Science, Chinese Academy of Sciences, Beijing 100190, PR China}

\date{}
\begin{abstract}
    We propose a unified Onsager-based neural time-stepping (ONTS) framework for 
    solving gradient-flow partial differential equations with conserved and nonconserved dynamics. 
    Unlike conventional neural time-stepping methods that directly parameterize 
    the evolving state, the proposed framework employs spatial neural networks 
    to approximate the Onsager process variables, namely, the rates in nonconserved systems and the fluxes in conserved systems. Within this Eulerian framework, 
    we develop three complementary schemes: an explicit least-squares scheme (E-ONTS), 
    an implicit variational scheme (I-ONTS), and a predictor--corrector scheme (PC-ONTS) 
    that combines explicit prediction with safeguarded implicit correction. For conserved 
    systems, the flux-based formulation guarantees discrete mass conservation under periodic or no-flux boundary conditions. Furthermore, we rigorously establish unconditional discrete energy stability for the implicit schemes, provided that a suitable descent criterion is satisfied, without requiring global minimization of the neural variational objective. Numerical experiments for the heat, Allen--Cahn, Cahn--Hilliard, 
    Fokker--Planck, and coupled Poisson--Nernst--Planck equations demonstrate 
    the accuracy, robustness, and structure-preserving properties of the proposed 
    framework. 
\end{abstract}


\end{frontmatter}
\section{Introduction}

In recent years, neural networks have emerged as flexible
computational tools for solving partial differential equations
(PDEs). Representative approaches include physics-informed
neural networks \cite{raissi2019physics, pang2019fpinns,
jagtap2020extended}, the deep Ritz and deep Galerkin
methods \cite{yu2018deep, liao2021, sirignano2018dgm},
weak-form neural network methods
\cite{zang2020weak, kharazmi2021hp}, deep learning
methods based on backward stochastic differential equations
\cite{han2017deep}, and random feature methods
\cite{chen2022bridging}. In parallel, operator learning
methods, such as DeepONet
\cite{lu2021learning, wang2021learning, wang2022improved},
Fourier neural operators \cite{li2021fourier}, and other
neural operator architectures \cite{kovachki2023neural,
gupta2021multiwavelet, li2024physics}, as well as deep
operator-splitting networks \cite{lan2023dosnet,
huang2027pidosnet}, have been developed to learn
solution operators of PDEs. These approaches have
demonstrated promising performance across a wide range
of problems.

However, approximation
accuracy alone does not necessarily guarantee that the
numerical solutions respect the physical and mathematical
structures of the underlying equations.
For time-dependent PDEs, properties such as conservation
laws and energy dissipation
are particularly important for capturing the correct
qualitative dynamics and long-time behavior of solutions.
Accordingly, increasing attention has been devoted to
incorporating such properties into neural-network-based
PDE solvers. For conservation laws, conservative physics-informed
neural networks (cPINNs) incorporate flux-continuity
conditions across subdomain interfaces
\cite{jagtap2020conservative}, while mass-preserving
spatio-temporal adaptive PINNs introduce mass conservation
constraints into the training objective for Cahn--Hilliard
equations \cite{huang2026mass}. 
Beyond PINN-based constraints, structure preservation has
also been pursued through compatible finite element
interpolated neural networks \cite{badia2025compatible}
and momentum-conserving neural solvers for elasticity
problems \cite{boon2025neural}.
To promote energy dissipation, energy-dissipation constraints
have been incorporated into PINNs for Allen--Cahn equations
\cite{kutuk2025energy}, while energy-stable neural networks
(EStable-Net) enforce discrete energy decay through
a block-based architecture \cite{wu2023energy}. 
At the operator-learning level, 
energy-dissipation-preserving
neural operators embed energy-stable schemes into network
architectures \cite{huang2026edp}, while the PENCO
framework incorporates physics, energy-dissipation, and
numerical-consistency constraints for three-dimensional
phase-field modeling \cite{bamdad2026penco}.
Other approaches include mass-conserving convolutional
neural networks based on fully discrete
operators \cite{geng2024deep}, energy-consistent neural
operators \cite{tanaka2025energy}, and the Sidecar framework,
which uses an auxiliary neural network to incorporate
conservation and dissipation properties
\cite{chen2026structure}.
Despite these developments, conventional PINNs typically
employ a unified space--time representation, with physical
constraints incorporated through additional loss terms.
Such a formulation makes it less straightforward to
leverage well-established structure-preserving time
discretizations with provable conservation, bound-preserving, 
or discrete energy-stability properties
\cite{eyre1998unconditionally,furihata2001stable,
shen2018scalar,hou2024energy, chen2025onsager, li2026high, zhang2026up, huang2026efficient}.

An alternative approach is to represent the solution using 
a spatial neural network whose parameters evolve in time, as in evolutional
deep neural networks (EDNNs)
\cite{du2021evolutional,kast2024positional}
and the closely related Neural Galerkin methods
\cite{bruna2024neural}.
These methods adopt a neural-network-based spatial
representation, $u(\bm{x},t)\approx u_{\bm{\theta}(t)}(\bm{x})$,
allowing spatial approximation and temporal integration
to be treated separately. Although this separation does not by itself guarantee 
the preservation of conservation laws or energy dissipation, it facilitates 
the incorporation of classical structure-preserving time discretizations. 
In this direction, related approaches have explored discrete
gradient-flow formulations for energy-stable time
integration \cite{georgoulis2023discrete}.
The separable deep minimizing movement (SDMM) method 
combines neural spatial representations with
variational time stepping for Allen--Cahn gradient flows
\cite{mattey2025gradient}.
Energy-dissipative evolutionary DeepONets
(EDE-DeepONets) employ the scalar auxiliary
variable (SAV) approach to ensure the unconditional
dissipation of a modified energy \cite{zhang2024energy}.
Within the broader Neural Galerkin setting, feature-based
methods preserve energy dissipation at the semi-discrete
level \cite{tang2026energy}, while nonlinear embeddings
and relaxation-projection techniques have been introduced
to preserve conserved quantities
\cite{schwerdtner2024nonlinear,shi2026preserving}.

A complementary perspective on structure preservation
is offered by Onsager's variational principle, which
characterizes dissipative dynamics in terms of free energy
and dissipation potentials.
This principle has also been extended to account for
inertial and thermal effects in non-isothermal two-phase
flows with phase transitions and wetting effects
\cite{qiao2026unified}.
For conserved systems, an Eulerian formulation can be
obtained by introducing flux variables subject to a
continuity equation \cite{chen2025onsager}.
Related variational neural approaches include energetic
variational neural networks (EVNNs), which employ an
Eulerian formulation for $L^2$-gradient flows but use
Lagrangian transport maps for conserved generalized
diffusions \cite{hu2024energetic}.
Onsager-type principles have also been used to learn
thermodynamically consistent dynamics, as well as free-energy
and dissipation potentials, from data
\cite{yu2021onsagernet,huang2022variational}.
More recently, Rayleighian minimization has been
employed to learn dissipative flux operators with
explicit time stepping \cite{chang2026unsupervised}.
Despite these advances, developing a unified Eulerian
neural time-stepping framework that accommodates both
conserved and nonconserved gradient flows, while preserving
mass conservation for conserved systems and discrete energy 
dissipation for both types of systems, remains an important challenge.

To address this challenge, we propose a unified Onsager-based neural 
time-stepping (ONTS) framework for conserved and nonconserved 
gradient-flow PDEs within an Eulerian formulation. Rather than 
directly parameterizing the evolving state, the proposed framework 
employs spatial neural networks to approximate the Onsager process 
variables, namely, the rates in nonconserved systems and fluxes in 
conserved systems, and advances the state through the corresponding 
kinematic relations. Within this framework, we develop three 
complementary schemes: an explicit least-squares Onsager 
neural time-stepping scheme (E-ONTS), an implicit Onsager 
neural time-stepping scheme (I-ONTS), and a predictor--corrector 
Onsager neural time-stepping scheme (PC-ONTS). The predictor--corrector 
scheme uses an explicit least-squares prediction to initialize the implicit 
variational correction, together with a safeguard for selecting 
the initial guess. For conserved systems, the flux-based formulation 
ensures mass conservation under periodic or no-flux boundary conditions. 
Moreover, we rigorously establish unconditional discrete energy stability 
for the implicit schemes under a suitable descent criterion, without 
requiring global minimization of the neural variational objective. Numerical 
experiments on heat, Allen--Cahn, Cahn--Hilliard, Fokker--Planck, 
and coupled Poisson--Nernst--Planck systems demonstrate the accuracy 
and structure-preserving properties of the proposed framework. 
In particular, the results confirm the predicted energy-dissipation 
behavior and demonstrate the robustness of the implicit schemes in 
stiff, high-dimensional, and coupled nonlinear problems.

The remainder of this paper is organized as follows.
Section 2 introduces the Onsager variational principle,
the associated time discretizations, and
neural-network-based spatial representations.
Section 3 presents the proposed Onsager neural
time-stepping framework, including the explicit,
implicit, and predictor--corrector schemes.
Section 4 evaluates the accuracy, efficiency, and
stability of the proposed schemes through numerical
experiments. Finally, Section 5 concludes the paper.

\section{Preliminaries}

In this section, we briefly review the main ingredients underlying the
proposed framework. We first introduce the Onsager variational principle (OVP) 
for nonconserved and conserved systems, followed by its time-discrete
formulations. We then review the neural-network-based spatial 
representation used in the proposed Onsager 
neural time-stepping schemes presented in Section~3.

\subsection{Onsager variational principle}

Consider a system described by the state variables 
$\bm{u}(\bm{x},t)=(u_1(\bm{x},t),\ldots,u_s(\bm{x},t))^T$ 
in $\Omega\subset\mathbb{R}^d$, with
$\partial_t\bm{u}=(\partial_t u_1,\ldots,\partial_t u_s)^T$
denoting its time derivative. 
For dissipative processes in the linear-response regime with
negligible inertial effects, the evolution equations can be
derived using the OVP \cite{chen2025onsager, doi2011onsager}.
We present the formulations for systems with nonconserved
and conserved variables separately below.

\subsubsection{Systems with nonconserved variables}

For a system described by nonconserved state variables,
let $\mathcal{E}(\bm{u})$ denote its free energy functional. 
Its rate of change along a sufficiently smooth trajectory 
$\bm{u}(t)$ is given by
\begin{equation*}
    \dot{\mathcal{E}}
    (\bm{u};\partial_t\bm{u})
    =
    \left\langle
        \frac{\delta \mathcal{E}(\bm{u})}
        {\delta \bm{u}},
        \partial_t\bm{u}
    \right\rangle .
\end{equation*}
Here, $\left\langle \cdot, \cdot \right\rangle$ 
denotes the $L^2(\Omega;\mathbb{R}^s)$ inner product, and 
$\delta\mathcal{E}/\delta\bm{u}$ denotes the variational 
derivative of the free energy. 
Under the linear-response assumption, the dissipation potential 
$\Phi(\bm{u};\cdot)$ is taken to be a positive-definite 
quadratic functional of the rate variable, consistent 
with Onsager's reciprocal relations. 
In particular, 
\begin{equation*}
    \Phi(\bm{u};\lambda\bm{v})
    =
    \lambda^2\Phi(\bm{u};\bm{v}),
    \qquad
    \Phi(\bm{u};\bm{v})\geq 0,
\end{equation*} 
for every $\lambda\in\mathbb{R}$ and admissible rate $\bm{v}$,
with equality in the second relation if and only if
$\bm{v}=\bm{0}$.

The Rayleighian is defined for an admissible rate $\bm{v}$ by
\begin{equation*}
    \mathcal{R}(\bm{u};\bm{v})
    :=
    \left\langle
        \frac{\delta\mathcal{E}(\bm{u})}{\delta\bm{u}},
        \bm{v}
    \right\rangle
    +
    \Phi(\bm{u};\bm{v}).
\end{equation*}
At each fixed state $\bm{u}$, the OVP determines the evolution
rate by minimizing the Rayleighian over all admissible rates:
\begin{equation*}
    \partial_t\bm{u}
    \in
    \arg\min_{\bm{v}}
    \mathcal{R}(\bm{u};\bm{v}).
\end{equation*}
The corresponding first-order optimality condition yields
\begin{equation*}
    \left.
    \frac{\delta \Phi(\bm{u};\bm{v})}
         {\delta \bm{v}}
    \right|_{\bm{v}=\partial_t\bm{u}}
    +
    \frac{\delta\mathcal{E}(\bm{u})}
         {\delta\bm{u}}
    =0,
\end{equation*}
which determines the evolution equation of the system. 
By the quadratic homogeneity of $\Phi$, this evolution satisfies
\begin{equation*}
    \frac{\mathrm{d}}{\mathrm{d}t}\mathcal{E}(\bm{u}) = -2 \Phi(\bm{u};\partial_t \bm{u}) \leq 0.
\end{equation*}
The Allen--Cahn equation is a standard example of this class of gradient
flows.

\subsubsection{Systems with conserved variables}

For a system described by conserved state variables,
each component $u_i$ satisfies a local conservation law
\begin{equation*}
    \partial_t u_i + \nabla\cdot \bm{j}_i = 0,
    \qquad i=1,\ldots,s,
\end{equation*}
where $\bm{j}_i$ denotes the flux associated with $u_i$.
Let $\bm{j} = (\bm{j}_1^T,\ldots,\bm{j}_s^T)^T$ 
denote the collection of fluxes. 
We assume periodic boundary conditions or, on a bounded
domain, the no-flux boundary conditions
$\bm{j}_i\cdot\bm{n}=0$ on $\partial\Omega$,
where $\bm{n}$ denotes the outward unit normal vector. 
Under the linear-response assumption, the dissipation potential 
$\Phi(\bm{u};\bm{j})$ is taken to be a positive-definite 
quadratic functional of the flux, consistent 
with Onsager's reciprocal relations, satisfying
\begin{equation*}
    \Phi(\bm{u};\lambda\bm{j})
    =
    \lambda^2
    \Phi(\bm{u};\bm{j}),
    \qquad
    \Phi(\bm{u};\bm{j})\geq 0,
\end{equation*}
for every $\lambda\in\mathbb{R}$ and admissible flux $\bm{j}$,
with equality in the second relation if and only if
$\bm{j}=\bm{0}$.

For a candidate rate $\bm{v}$ and flux $\bm{q}$, the
Rayleighian is defined as
\begin{equation*}
    \mathcal{R}(\bm{u};\bm{v},\bm{q})
    :=
    \left\langle
        \frac{\delta\mathcal{E}(\bm{u})}{\delta\bm{u}},
        \bm{v}
    \right\rangle
    +
    \Phi(\bm{u};\bm{q}).
\end{equation*}
At each fixed state $\bm{u}$, the OVP determines the evolution
rate and flux by minimizing the Rayleighian subject to the
conservation laws and the no-flux boundary conditions:
\begin{equation*}
    (\partial_t\bm{u},\bm{j})
    \in
    \arg\min_{\bm{v},\bm{q}}
    \left\{
        \mathcal{R}
        (\bm{u};\bm{v},\bm{q})
        :
        \bm{v}
        +
        \nabla\cdot\bm{q}
        =0
    \right\}.
\end{equation*}
Here, the minimization is restricted to fluxes satisfying
the prescribed periodic or no-flux boundary conditions,
and $\nabla\cdot\bm{q}=(\nabla\cdot\bm{q}_1,\ldots,\nabla\cdot\bm{q}_s)^T$.

Introducing a vector-valued Lagrange multiplier field
$\bm{\mu}=(\mu_1,\ldots,\mu_s)^T$, we define the Lagrangian
\begin{equation*}
    \mathcal{R}_{\bm{\mu}}(\bm{u};\bm{v},\bm{q})
    :=
    \mathcal{R}(\bm{u};\bm{v},\bm{q})
    -
    \sum_{i=1}^{s}
    \int_{\Omega}
    \mu_i
    \left(
        v_i+\nabla\cdot\bm{q}_i
    \right)
    \,\mathrm{d}\bm{x}.
\end{equation*}
Taking variations with respect to $\bm{v}$, $\bm{q}$, and
$\bm{\mu}$, and integrating by parts using the no-flux
boundary conditions, yields the following first-order
optimality conditions at
$(\bm{v},\bm{q})=(\partial_t\bm{u},\bm{j})$:
\begin{equation*}
    \begin{cases}
        \partial_t u_i
        +
        \nabla\cdot\bm{j}_i
        =0,
        \\[2mm]
        \displaystyle
        \frac{\delta\Phi}
        {\delta\bm{j}_i}
        +
        \nabla\mu_i
        =0,
        \\[3mm]
        \displaystyle
        \frac{\delta\mathcal{E}}
        {\delta u_i}
        -
        \mu_i
        =0,
    \end{cases}
    \qquad i=1,\ldots,s.
\end{equation*} 
Thus, the Lagrange multipliers are identified with the
chemical potentials, and the second equation determines
the constitutive relation for the fluxes. Combining these 
relations with the conservation laws and using the quadratic homogeneity
of $\Phi$, we obtain
\begin{equation*}
    \frac{\mathrm{d}}{\mathrm{d}t}\mathcal{E}(\bm{u}) = -2 \Phi(\bm{u};\bm{j}) \leq 0.
\end{equation*}
Representative examples include 
the heat, Cahn--Hilliard, Fokker--Planck, 
and Poisson--Nernst--Planck equations.

\subsection{Time-discretization based on the Onsager principle}

We next recall the time-discrete formulations based on the OVP
\cite{chen2025onsager}. For simplicity, we consider a uniform partition
$t_k=k\tau$ of $[0,T]$, where $\tau=T/N>0$ and $k=0,\ldots,N$.
Let $\bm{u}^k$ denote the approximation at time $t_k$.
All minimizations below are understood to be over admissible states and
fluxes for which the energy and dissipation potentials are well defined;
in particular, the previous state is assumed to be admissible.
In both formulations, the state argument of the dissipation potential is
evaluated at $\bm{u}^k$.

\subsubsection{Discretization for systems with nonconserved variables}

Given $\bm{u}^k$, we define the discrete Rayleighian
for a candidate state $\bm{u}$ by
\begin{equation*}
    \mathcal{R}_\tau^k(\bm{u}^k;\bm{u})
    :=
    D_\tau\mathcal{E}(\bm{u}^k;\bm{u})
    +
    \Phi\left(
        \bm{u}^k;
        \frac{\bm{u}-\bm{u}^k}{\tau}
    \right),
\end{equation*}
where
\begin{equation*}
    D_\tau\mathcal{E}(\bm{u}^k;\bm{u})
    :=
    \frac{\mathcal{E}(\bm{u})-\mathcal{E}(\bm{u}^k)}{\tau}.
\end{equation*}
The time-discrete OVP defines the next state $\bm{u}^{k+1}$ through
\begin{equation*}
    \bm{u}^{k+1}\in\arg\min_{\bm{u}}\mathcal{R}_\tau^k(\bm{u}^k;\bm{u}).
\end{equation*}
By the quadratic homogeneity of $\Phi$, this is equivalent to
\begin{equation*}
    \bm{u}^{k+1}
    \in
    \arg\min_{\bm{u}}
    \left\{
        \mathcal{E}(\bm{u})
        +
        \frac{1}{\tau}
        \Phi(\bm{u}^k;\bm{u}-\bm{u}^k)
    \right\}.
\end{equation*}
This is a minimizing-movement-type scheme with a quadratic dissipation
penalty evaluated at the previous state. 
If a minimizer exists, comparison with the admissible candidate $\bm{u}^k$ yields
\begin{equation*}
    \mathcal{E}(\bm{u}^{k+1})
    +
    \frac{1}{\tau}
    \Phi(\bm{u}^k;\bm{u}^{k+1}-\bm{u}^k)
    \leq
    \mathcal{E}(\bm{u}^k),
\end{equation*}
which guarantees energy dissipation for any
time-step size $\tau>0$.

\subsubsection{Discretization for systems with conserved variables}

For conserved systems, the continuity equations are discretized as
\begin{equation*}
    \frac{u_i-u_i^k}{\tau}+\nabla\cdot \bm{j}_i=0,
    \qquad i=1,\ldots,s,
\end{equation*}
where $\bm{j}_i$ represents a piecewise-constant-in-time
approximation of the flux over $[t_k,t_{k+1}]$
and satisfies the prescribed periodic or no-flux
boundary conditions. The discrete Rayleighian is
\begin{equation*}
    \mathcal{R}_\tau^k(\bm{u}^k;\bm{u},\bm{j})
    :=
    D_\tau\mathcal{E}(\bm{u}^k;\bm{u})
    +
    \Phi(\bm{u}^k;\bm{j}).
\end{equation*}
The discrete OVP determines the next state and flux by
\begin{equation*}
    (\bm{u}^{k+1},\bm{j}^{k+1})
    \in
    \arg\min_{\bm{u},\bm{j}}
    \left\{
        \mathcal{R}_\tau^k(\bm{u}^k;\bm{u},\bm{j})
        :
        \frac{\bm{u}-\bm{u}^k}{\tau}
        +\nabla\cdot\bm{j}=0
    \right\}.
\end{equation*}
Here, the minimization is restricted to fluxes satisfying
the prescribed periodic or no-flux boundary conditions. 
Equivalently,
\begin{equation*}
    (\bm{u}^{k+1},\bm{j}^{k+1})
    \in
    \arg\min_{\bm{u},\bm{j}}
    \left\{
        \mathcal{E}(\bm{u})
        +
        \tau\Phi(\bm{u}^k;\bm{j})
        :
        \frac{\bm{u}-\bm{u}^k}{\tau}
        +\nabla\cdot\bm{j}=0
    \right\}.
\end{equation*}

Comparison with the admissible pair $(\bm{u}^k,\bm{0})$ yields
\begin{equation*}
    \mathcal{E}(\bm{u}^{k+1})
    +\tau\Phi(\bm{u}^k;\bm{j}^{k+1})
    \leq\mathcal{E}(\bm{u}^k).
\end{equation*}
Moreover, integrating the discrete continuity equations
over $\Omega$ and using the periodic or no-flux
boundary conditions gives 
\begin{equation*}
    \int_{\Omega}u_i^{k+1}\,\mathrm{d}\bm{x}
    =\int_{\Omega}u_i^k\,\mathrm{d}\bm{x},
    \qquad i=1,\ldots,s.
\end{equation*}
Thus, conservation follows from the discrete continuity equations,
whereas energy dissipation follows from the variational comparison
principle \cite{chen2025onsager}. Importantly, in both 
formulations, the corresponding energy-dissipation 
inequality remains valid for any admissible approximate 
minimizer whose objective value does not exceed that of 
the unchanged-state candidate. Thus, exact minimization 
is not required to ensure discrete energy stability.

These time-discrete variational formulations provide the basis for
structure-preserving schemes. In the following subsection, we recall the
neural-network-based spatial approximation used to implement the resulting
time-stepping methods.

\subsection{Neural-network-based spatial representation}

Consider a time-dependent PDE of the form
\begin{equation}
    \partial_t\bm{u}=\mathcal{N}_{\bm{x}}(\bm{u}),
    \qquad \bm{x}\in\Omega\subset\mathbb{R}^d,
    \label{eq:general_pde}
\end{equation}
where $\bm{u}(\bm{x},t)\in\mathbb{R}^s$ and
$\mathcal{N}_{\bm{x}}$ is a possibly nonlinear spatial differential
operator. In evolutional neural-network approaches, the spatial dependence
of the solution is represented by a neural network with time-dependent
parameters \cite{du2021evolutional,kast2024positional,bruna2024neural}:
\begin{equation*}
    \bm{u}(\bm{x},t)
    \approx
    \hat{\bm{u}}
    \bigl(\bm{x};\bm{\theta}(t)\bigr),
\end{equation*}
where $\bm{\theta}(t)\in\mathbb{R}^p$ collects the trainable parameters. 
Therefore, the neural network acts as a spatial representation of the
solution, while its temporal dependence is encoded through the evolution
of $\bm{\theta}(t)$. 
By the chain rule, the time derivative of the neural representation is
given by
\begin{equation*}
    \partial_t
    \hat{\bm{u}}
    \bigl(\bm{x};\bm{\theta}(t)\bigr)
    =
    D_{\bm{\theta}}
    \hat{\bm{u}}
    \bigl(\bm{x};\bm{\theta}(t)\bigr)
    \dot{\bm{\theta}}(t).
\end{equation*}

In the standard EDNNs formulation, the parameter velocity
$\dot{\bm{\theta}}$ is determined by minimizing the instantaneous PDE residual. 
At fixed $\bm{\theta}(t)$,
\begin{equation*}
    \dot{\bm{\theta}}
    \in
    \arg\min_{\bm{\eta}\in\mathbb{R}^{p}}
    \frac{1}{2}
    \int_{\Omega}
    \left\|
        D_{\bm{\theta}}
        \hat{\bm{u}}
        \,\bm{\eta}
        -
        \mathcal{N}_{\bm{x}}
        (\hat{\bm{u}})
    \right\|_2^2
    \,\mathrm{d}\bm{x}.
\end{equation*}
The corresponding first-order optimality condition is
\begin{equation*}
    \left(
        \int_{\Omega}
        \left(D_{\bm{\theta}}\hat{\bm{u}}\right)^T
        D_{\bm{\theta}}\hat{\bm{u}}
        \,\mathrm{d}\bm{x}
    \right)
    \dot{\bm{\theta}}
    =
    \int_{\Omega}
    \left(D_{\bm{\theta}}\hat{\bm{u}}\right)^T
    \mathcal{N}_{\bm{x}}(\hat{\bm{u}})
    \,\mathrm{d}\bm{x},
\end{equation*}
where $D_{\bm{\theta}}\hat{\bm{u}}\in\mathbb{R}^{s\times p}$
denotes the Jacobian of the neural representation with respect to its
parameters. The resulting Gram matrix may be singular or ill-conditioned
if the neural-network parameterization is redundant.

In practice, the spatial integrals are approximated by quadrature or
Monte Carlo sampling. To illustrate the equal-weight case, let
$\{\bm{x}_i\}_{i=1}^{N_c}$ be the collocation points and define
\begin{equation*}
    \bm{J}
    =
    \begin{bmatrix}
        D_{\bm{\theta}}
        \hat{\bm{u}}(\bm{x}_1;\bm{\theta})
        \\
        \vdots
        \\
        D_{\bm{\theta}}
        \hat{\bm{u}}(\bm{x}_{N_c};\bm{\theta})
    \end{bmatrix},
    \qquad
    \bm{N}
    =
    \begin{bmatrix}
        \mathcal{N}_{\bm{x}}
        (\hat{\bm{u}})(\bm{x}_1)
        \\
        \vdots
        \\
        \mathcal{N}_{\bm{x}}
        (\hat{\bm{u}})(\bm{x}_{N_c})
    \end{bmatrix}.
\end{equation*}
Omitting the common positive quadrature factor, which does not affect the
minimizer, the discrete least-squares problem is
\begin{equation*}
    \dot{\bm{\theta}}
    \in
    \arg\min_{\bm{\eta}\in\mathbb{R}^{p}}
    \frac{1}{2}
    \left\|
        \bm{J}\bm{\eta}
        -
        \bm{N}
    \right\|_2^2,
\end{equation*}
with normal equations
\begin{equation*}
    \bm{J}^{T}\bm{J}\,
    \dot{\bm{\theta}}
    =
    \bm{J}^{T}\bm{N}.
\end{equation*}
The parameter derivatives in $\bm{J}$ and the spatial derivatives
required in $\bm{N}$ can be computed by automatic differentiation. 
Once the parameter velocity $\dot{\bm{\theta}}$ is 
determined, the network parameters are advanced 
using a numerical ODE integrator. For example, forward Euler
gives
\begin{equation*}
    \bm{\theta}^{k+1}
    =
    \bm{\theta}^{k}
    +
    \tau\dot{\bm{\theta}}^{k},
\end{equation*}
while higher-order schemes, such as Runge--Kutta methods, may also be
used. The initial parameters $\bm{\theta}^{0}$ are typically obtained by
fitting the prescribed initial condition.

The standard EDNNs formulation evolves the parameters of a
neural-network representation of the state through a
PDE-residual projection. In contrast, the proposed ONTS
framework employs spatial neural networks to approximate
the Onsager process variables, namely rates for nonconserved
systems and fluxes for conserved systems, and advances the
state through the corresponding kinematic relations.
The resulting explicit, implicit, and predictor--corrector
Onsager neural time-stepping schemes are developed in
Section~3.

\section{A unified Onsager neural time-discretization framework}

In this section, we combine the time-discrete Onsager
variational formulations with neural-network-based spatial
representations to develop a unified numerical framework
for dissipative systems. Unlike conventional EDNNs methods \cite{du2021evolutional},
which directly parameterize the state variables, our approach
employs neural networks to represent the Onsager process
variables, while the evolution of the physical state is
determined by the corresponding kinematic relations.

We develop two basic schemes and a predictor--corrector
strategy. The explicit least-squares Onsager neural
time-stepping (E-ONTS) scheme advances the state using
the current neural process variable, constructs an
Onsager-consistent target at the updated state, and
updates the network parameters through a linearized
least-squares fit. The implicit Onsager neural time-stepping
(I-ONTS) scheme determines the neural process variable
by minimizing the time-discrete Onsager functional.
The predictor--corrector Onsager neural time-stepping
(PC-ONTS) scheme combines explicit prediction with
implicit variational correction, using an additional
Rayleighian-based safeguard for multidimensional problems.

\subsection{Neural spatial representation}
\label{Neural_spatial_representation}
We employ neural networks to represent the spatial
dependence of the Onsager process variables. The
specific quantity represented by the network depends
on whether the state variables are nonconserved or
conserved. At each time step, let
$\bm{\theta}\in\mathbb{R}^{p}$ denote the trainable
network parameters. 

For systems with nonconserved variables, the Onsager
process variable is the temporal rate
$\bm{v}=\partial_t\bm{u}$.
We introduce the neural representation
\begin{equation*}
    \bm{v}_{\bm{\theta}(t)}(\bm{x})\approx \bm{v}(\bm{x},t) = \partial_t \bm{u},
\end{equation*}
where the neural parameters are determined at each 
time step. Given the state $\bm{u}^{k}$ at time $t_k$, a candidate
neural rate $\bm{v}_{\bm{\theta}}$ defines the trial state
at the next time level by
\begin{equation*}
    \bm{u}^{k+1}_{\tau}
    (\bm{\theta})
    =
    \bm{u}^{k}
    +
    \tau
    \bm{v}_{\bm{\theta}},
\end{equation*}
where $\tau=t_{k+1}-t_k>0$.

For systems with conserved variables, the evolution
is subject to the continuity equation
\begin{equation*}
    \partial_t\bm{u}
    +
    \nabla\cdot\bm{j}
    =
    0,
\end{equation*}
where $\bm{j}$ collects the fluxes associated with the
components of $\bm{u}$. In this case,
the Onsager process variable is
the physical flux, which is approximated by
\begin{equation*}
    \bm{j}_{\bm{\theta}(t)}(\bm{x})
    \approx \bm{j}(\bm{x},t).
\end{equation*}
The corresponding trial state is determined through
the discrete continuity equation,
\begin{equation*}
    \bm{u}^{k+1}_{\tau}(\bm{\theta})
    =
    \bm{u}^{k}
    -
    \tau\nabla\cdot\bm{j}_{\bm{\theta}},
\end{equation*}
where the divergence is applied componentwise
to the fluxes.

Both nonconserved and conserved systems can be
described using a unified neural process variable,
\begin{equation}
    \bm{z}_{\bm{\theta}}
    =
    \begin{cases}
        \bm{v}_{\bm{\theta}},
        & \text{for nonconserved variables},\\[1mm]
        \bm{j}_{\bm{\theta}},
        & \text{for conserved variables}.
    \end{cases}
    \label{eq:unified-neural-rate}
\end{equation}
Correspondingly, the trial-state map takes
the unified form
\begin{equation}
    \bm{u}_{\tau}^{k+1}(\bm{\theta})
    =
    \bm{u}^{k}
    +
    \tau\mathcal{B}\bm{z}_{\bm{\theta}},
    \label{eq:unified-neural-update}
\end{equation}
where the kinematic operator $\mathcal{B}$ is defined by
\begin{equation}
    \mathcal{B}\bm{z}
    =
    \begin{cases}
        \bm{v},
        & \text{for nonconserved variables},\\[1mm]
        -\nabla\cdot\bm{j},
        & \text{for conserved variables}.
    \end{cases}
    \label{eq:kinematic-operator}
\end{equation}
The neural process variables are assumed to satisfy
the relevant kinematic boundary conditions, and
the associated trial states are required to belong
to the admissible state space. In particular, we
assume that the neural representation admits the
zero process variable, which provides a feasible
reference state for the implicit variational schemes.

This unified representation allows the physical
state to be updated through its kinematic relation,
while the neural-network parameters determine
the corresponding process variables. It serves
as the common spatial approximation framework
for the explicit, implicit, and predictor--corrector
Onsager time-stepping schemes developed below.

\subsection{Explicit least-squares Onsager neural time-stepping (E-ONTS)}

We first construct an explicit realization of the Onsager neural
framework using the process variable representation introduced in the
previous subsection. At each time step, the current neural process
variable is used to advance the physical state. An Onsager-consistent
target is then evaluated at this updated state, and the network
parameters are updated by a linearized least-squares fit to the target.
The least-squares step does not minimize the time-discrete Onsager
functional; rather, the Onsager principle specifies the target process
variable, while the least-squares projection updates its neural
representation. The initial parameters $\bm\theta^0$ are obtained by
fitting the neural process variable to an Onsager-consistent target
evaluated at the prescribed initial state $\bm u^0$.

\paragraph{Nonconserved variables}
For nonconserved variables, the neural process variable is the rate
$\bm v_{\bm\theta}\approx\partial_t\bm u$. Given $\bm u^k$ and
$\bm\theta^k$, we first update the state explicitly
\begin{equation}
    \bm{u}^{k+1}_{\tau}
    =
    \bm{u}^{k}
    +
    \tau
    \bm{v}_{\bm{\theta}^k}.
    \label{eq:nonconserved-update-state}
\end{equation}
where $\tau=t_{k+1}-t_k>0$. At the updated state $\bm{u}^{k+1}_{\tau}$, 
the Onsager-consistent target rate is defined by
\begin{equation*}
    \bm v_{\mathrm{tar}}^{k+1}
    \in\arg\min_{\bm v\in\mathcal V_{\mathrm{ad}}}
    \left\{
        \left\langle
            \frac{\delta\mathcal E}{\delta\bm{u}}(\bm{u}^{k+1}_{\tau}),
            \bm v
        \right\rangle
        +\Phi(\bm{u}^{k+1}_{\tau};\bm v)
    \right\},
\end{equation*}
where $\mathcal V_{\mathrm{ad}}$ denotes the admissible rate space.
Assuming differentiability and an unconstrained interior minimizer,
the corresponding first-order optimality condition is
\begin{equation*}
    \frac{\delta\mathcal E}{\delta\bm{u}}(\bm{u}^{k+1}_{\tau})
    +\frac{\delta\Phi}{\delta\bm v}
       (\bm{u}^{k+1}_{\tau};\bm v_{\mathrm{tar}}^{k+1})=\bm0.
\end{equation*}

To fit the target $\bm v_{\mathrm{tar}}^{k+1}$,
we linearize the neural rate with respect to its
parameters around $\bm\theta^k$:
\begin{equation*}
    \bm v_{\bm\theta^k+\Delta\bm\theta}
    \approx
    \bm v_{\bm\theta^k}
    +
    D_{\bm\theta}\bm v_{\bm\theta^k}
    \,\Delta\bm\theta.
\end{equation*}
Writing $\Delta\bm\theta=\tau\bm\eta$,
we determine the parameter velocity from
the linear least-squares problem
\begin{equation}
    \dot{\bm{\theta}}^k
    \in
    \arg\min_{\bm{\eta}\in\mathbb R^p}
    \left\|
        D_{\bm{\theta}}
        \bm v_{\bm{\theta}^k}
        \bm\eta
        -
        \frac{
            \bm v_{\mathrm{tar}}^{k+1}
            -
            \bm v_{\bm{\theta}^k}
        }{\tau}
    \right\|_{L^2(\Omega)}^2 .
    \label{eq:nonconserved-explicit-ls}
\end{equation}
The spatial norm is approximated by quadrature or sampling, giving a
matrix least-squares problem. Finally, the parameters are advanced by
\begin{equation*}
    \bm{\theta}^{k+1}
    =
    \bm{\theta}^{k}
    +
    \tau
    \dot{\bm{\theta}}^k .
\end{equation*}
Notice that $\bm{u}^{k+1}_{\tau}$ in
\eqref{eq:nonconserved-update-state} is computed using
$\bm\theta^k$; the updated network represents the rate for the
next explicit state update and does not retroactively alter
$\bm{u}^{k+1}_{\tau}$.

\paragraph{Conserved variables}
For conserved variables, the process variable is the physical flux,
which satisfies the continuity equation
\begin{equation*}
    \partial_t\bm u+\nabla\cdot\bm j=\bm0.
\end{equation*}
We use the neural approximation $\bm j_{\bm\theta}$ and assume
that it satisfies the prescribed periodic or no-flux boundary
conditions for all admissible parameter values. Given $\bm u^k$ and
$\bm\theta^k$, the explicit conservative update is
\begin{equation}
    \bm u^{k+1}_{\tau}=\bm u^k-\tau\nabla\cdot\bm j_{\bm\theta^k}.
    \label{eq:explicit-conserved-state-update}
\end{equation}
The target flux at the updated state is determined by the Onsager
Rayleighian:
\begin{equation*}
    \bm j_{\mathrm{tar}}^{k+1}
    \in\arg\min_{\bm j\in\mathcal J_{\mathrm{ad}}}
    \left\{
       \left\langle
          \frac{\delta\mathcal E}{\delta\bm u}(\bm u^{k+1}_{\tau}),
          -\nabla\cdot\bm j
       \right\rangle
       +\Phi(\bm u^{k+1}_{\tau};\bm j)
    \right\},
\end{equation*}
where $\mathcal J_{\mathrm{ad}}$ contains fluxes satisfying the
required boundary conditions. Under periodic or no-flux conditions,
integration by parts rewrites the energy-rate term as
\begin{equation*}
    \left\langle
       \nabla\frac{\delta\mathcal E}{\delta\bm u}(\bm u^{k+1}_{\tau}),
       \bm j
    \right\rangle.
\end{equation*}
Here and below, gradients and divergences act componentwise on the
corresponding state and flux components.

As in the nonconserved case, we linearize the neural flux about
$\bm\theta^k$ and fit the target in the resulting affine space:
\begin{equation}
    \dot{\bm\theta}^{k}\in
    \arg\min_{\bm\eta\in\mathbb R^p}
    \left\|
       D_{\bm\theta}\bm j_{\bm\theta^k}\bm\eta
       -\frac{\bm j_{\mathrm{tar}}^{k+1}-\bm j_{\bm\theta^k}}{\tau}
    \right\|_{L^2(\Omega)}^2.
    \label{eq:conserved-explicit-ls}
\end{equation}
We then update
\begin{equation*}
    \bm\theta^{k+1}=\bm\theta^k+\tau\dot{\bm\theta}^{k}.
\end{equation*}
The resulting flux $\bm j_{\bm\theta^{k+1}}$ is used in the next
state update, from $\bm{u}^{k+1}$ to $\bm{u}^{k+2}$.

\begin{remark}[Alternative Rayleighian-based projection]
The unweighted $L^2$ target fits in
\eqref{eq:nonconserved-explicit-ls} and
\eqref{eq:conserved-explicit-ls} can be replaced by a projection
that minimizes the continuous Rayleighian over an affine
linearization of the neural process variable representation. 
When an explicit constitutive relation is unavailable,
this approach avoids computing the full-space target
rate or flux before the parameter update. 
For nonconserved variables, this gives
\begin{equation*}
    \dot{\bm\theta}^{k}\in\arg\min_{\bm\eta\in\mathbb R^p}
    \mathcal R\left(
       \bm u^{k+1}_{\tau};
       \bm v_{\bm\theta^k}
       +\tau D_{\bm\theta}\bm v_{\bm\theta^k}\bm\eta
    \right),
\end{equation*}
whereas, for conserved variables, it gives
\begin{equation*}
    \dot{\bm\theta}^{k}\in\arg\min_{\bm\eta\in\mathbb R^p}
    \Bigl\{\left\langle
       \frac{\delta\mathcal E}{\delta\bm u}(\bm u^{k+1}_{\tau}),
       -\nabla\cdot\bm j_{\mathrm{aff}}(\bm\eta)
    \right\rangle+\Phi\bigl(\bm u^{k+1}_{\tau};\bm j_{\mathrm{aff}}(\bm\eta)\bigr)
    \Bigr\},
\end{equation*}
with
$\bm j_{\mathrm{aff}}(\bm\eta)
 :=\bm j_{\bm\theta^k}
   +\tau D_{\bm\theta}\bm j_{\bm\theta^k}\bm\eta$.
The affine fluxes must satisfy the prescribed boundary conditions.
For a positive-definite quadratic dissipation potential, these
minimizations are equivalent to least-squares projections in the
metric induced by the dissipation operator, which generally differ
from the unweighted $L^2$ fits. They are alternative updates, not
additional steps in the E-ONTS algorithm below.
\end{remark}

\begin{proposition}[Mass conservation of E-ONTS]
\label{prop:explicit-mass-conservation}
Assume that each component of the neural flux is sufficiently regular
and satisfies either a no-flux condition
$\bm j_{i,\bm\theta^k}\cdot\bm n=0$ on $\partial\Omega$
or the corresponding periodic boundary conditions. Then the state
update \eqref{eq:explicit-conserved-state-update} satisfies
\begin{equation*}
    \int_\Omega u_i^{k+1}\,\mathrm{d}\bm x
    =\int_\Omega u_i^k\,\mathrm{d}\bm x,
    \qquad i=1,\ldots,s,
\end{equation*}
for every $\tau>0$.
\end{proposition}

\begin{proof}
Integrating \eqref{eq:explicit-conserved-state-update} over
$\Omega$ and applying the divergence theorem yields
\begin{align*}
    \int_\Omega u_i^{k+1}\,\mathrm{d}\bm x
    &=\int_\Omega u_i^k\,\mathrm{d}\bm x
      -\tau\int_\Omega\nabla\cdot
            \bm j_{i,\bm\theta^k}\,\mathrm{d}\bm x\\
    &=\int_\Omega u_i^k\,\mathrm{d}\bm x
      -\tau\int_{\partial\Omega}
            \bm j_{i,\bm\theta^k}\cdot\bm n\,\mathrm{d}S\\
    &=\int_\Omega u_i^k\,\mathrm{d}\bm x.
\end{align*}
For periodic boundary conditions, opposite boundary fluxes cancel.
\end{proof}

Thus, E-ONTS inherits mass conservation from the conservative state
update, provided that the flux representation satisfies the stated
boundary conditions. Using the unified process variable $\bm z_{\bm\theta}$ and kinematic
operator $\mathcal B$ from
\eqref{eq:unified-neural-rate}--\eqref{eq:kinematic-operator},
Algorithm~\ref{alg:e-onts} summarizes the explicit scheme for both
nonconserved and conserved systems. The target minimization is understood
over the appropriate admissible process-variable space and can be
replaced by an explicit constitutive expression when available.
E-ONTS avoids minimizing the nonlinear
time-discrete Onsager functional, but evaluating the target and solving
the linearized least-squares problem still incur computational cost.
The explicit state update does not, in general, guarantee discrete
energy dissipation and may require a step size restriction. 
This motivates the implicit
Onsager neural time-stepping (I-ONTS) scheme introduced next.

\begin{algorithm}[htbp]
\caption{Explicit least-squares Onsager neural time-stepping (E-ONTS)}
\label{alg:e-onts}
\hspace*{0.02in} {\bf Require:}
        Initial state $\bm u^0$; step size $\tau>0$; number of steps $N$;
         process network $\bm z_{\bm\theta}$; kinematic operator $\mathcal B$;
         free energy $\mathcal E$ and dissipation potential $\Phi$.
\begin{algorithmic}[1]
\State Compute an Onsager-consistent target $\bm z_{\mathrm{tar}}^0$ at $\bm u^0$.
\State Fit $\bm\theta^0$ via neural network training so that $\bm z_{\bm\theta^0}\approx\bm z_{\mathrm{tar}}^0$.
\For{$k=0,1,\ldots,N-1$}
    \State $\bm u^{k+1}_{\tau}\gets\bm u^k+\tau\mathcal B\bm z_{\bm\theta^k}$, $\bm\mu^{k+1}_{\tau}\gets\delta\mathcal E(\bm u^{k+1}_{\tau})/\delta\bm u$
    \State Determine the Onsager target 
    \begin{equation*}
        \displaystyle
      \bm z_{\mathrm{tar}}^{k+1}\in
      \arg\min_{\bm z\in\mathcal Z_{\mathrm{ad}}}
      \left\{\langle\bm\mu^{k+1}_{\tau},\mathcal B\bm z\rangle
      +\Phi(\bm u^{k+1}_{\tau};\bm z)\right\}
    \end{equation*}

    \State Let $\bm J^k= D_{\bm\theta}\bm z_{\bm\theta^k}$ and determine the parameter velocity
    \begin{equation*}
      \displaystyle
      \bm\eta^k\in\arg\min_{\bm\eta\in\mathbb R^p}
      \left\|\bm J^k\bm\eta
      -\frac{\bm z_{\mathrm{tar}}^{k+1}-\bm z_{\bm\theta^k}}{\tau}
      \right\|_{L^2(\Omega)}^2
    \end{equation*}
    \State $\bm\theta^{k+1}\gets\bm\theta^k+\tau\bm\eta^k$
\EndFor
\end{algorithmic}
\hspace*{0.02in} {\bf Output:} States $\{\bm u^k\}_{k=0}^{N}$ and process parameters $\{\bm\theta^k\}_{k=0}^{N}$
\end{algorithm}

\subsection{Implicit Onsager neural time-stepping (I-ONTS)}

The implicit scheme retains the time-discrete Onsager variational
structure in the neural spatial representation. Unlike E-ONTS, it
determines the process variable directly through a nonlinear
variational minimization, without separately evaluating an
Onsager-consistent target and evolving the parameters through a
linearized least-squares system.

For a given state $\bm u^k$, let $\bm z_{\bm\theta}$
and $\mathcal B$ denote the unified neural process
variable and kinematic operator introduced in
Section~\ref{Neural_spatial_representation}. 
Following the time-discrete Onsager
variational principle, we define the stepwise objective
\begin{equation}
    \mathcal J_k(\bm\theta)
    :=
    \mathcal E\bigl(
        \bm u_\tau^{k+1}(\bm\theta)
    \bigr)
    +
    \tau\Phi(\bm u^k;\bm z_{\bm\theta}),
    \label{eq:implicit-unified-objective}
\end{equation}
where
$\bm u_\tau^{k+1}(\bm\theta)
=\bm u^k+\tau\mathcal B\bm z_{\bm\theta}$. 
Throughout, the admissible parameters are those
for which the process variable satisfies the
prescribed boundary conditions and the trial state
belongs to the domain of the free energy.
We assume that the current state is admissible
and that the neural representation admits the zero
process variable, i.e.,
$\bm z_{\bm\theta_0}=\bm0$
for some $\bm\theta_0$.
By the quadratic structure of the dissipation potential,
$\Phi(\bm u^k;\bm0)=0$, and hence
$\mathcal J_k(\bm\theta_0)=\mathcal E(\bm u^k)$.
In practice, the minimization is performed approximately
using a gradient-based optimizer. 

\paragraph{Nonconserved variables}
For nonconserved variables, the network $\bm v_{\bm\theta}$ represents 
the temporal rate $\partial_t \bm{u}$. Its trial state is
\begin{equation*}
    \bm u^{k+1}_{\tau}(\bm\theta)
    =\bm u^k+\tau\bm v_{\bm\theta}
\end{equation*}
and the ideal implicit update is defined by
\begin{equation}
    \bm\theta^{k+1}\in\arg\min_{\bm\theta\in\Theta_{\mathrm{ad}}^k}
    \left\{
       \mathcal E\bigl(\bm u^{k+1}_{\tau}(\bm\theta)\bigr)
       +\tau\Phi(\bm u^k;\bm v_{\bm\theta})
    \right\},
    \label{eq:implicit-nonconserved-minimization}
\end{equation}
where $\Theta_{\mathrm{ad}}^k$ denotes the admissible parameter set
at step $k$. Once a parameter has been selected, the state is updated by
\begin{equation}
    \bm u^{k+1}=\bm u^k+\tau\bm v_{\bm\theta^{k+1}}.
    \label{eq:implicit-nonconserved-update}
\end{equation}
Unlike the explicit scheme, the trial state and the neural rate are
coupled within the same variational objective.

\begin{proposition}[Energy dissipation for nonconserved variables]
\label{prop:implicit-nonconserved-energy}
Let $\tau>0$, and suppose that the zero neural rate is admissible,
$\bm v_{\bm\theta_0}=\bm0$. Assume that
$\Phi(\bm u;\bm0)=0$ and $\Phi(\bm u;\bm v)\geq0$ for every
admissible rate. If an admissible parameter $\bm\theta^{k+1}$,
not necessarily a global minimizer of
\eqref{eq:implicit-nonconserved-minimization},
satisfies
\begin{equation}
    \mathcal J_k(\bm\theta^{k+1})
    \leq
    \mathcal J_k(\bm\theta_0)
    =
    \mathcal E(\bm u^k),
    \label{eq:implicit-nonconserved-descent}
\end{equation}
then the update
\eqref{eq:implicit-nonconserved-update}
obeys
\begin{equation}
    \mathcal E(\bm u^{k+1})
    +\tau\Phi(\bm u^k;\bm v_{\bm\theta^{k+1}})
    \leq\mathcal E(\bm u^k),
    \label{eq:implicit-nonconserved-energy-law}
\end{equation}
and, in particular,
\begin{equation}
    \mathcal E(\bm u^{k+1})\leq\mathcal E(\bm u^k).
    \label{eq:implicit-nonconserved-energy-decay}
\end{equation}
\end{proposition}

\begin{proof}
Substitution of \eqref{eq:implicit-nonconserved-update} into
\eqref{eq:implicit-nonconserved-descent} gives
\eqref{eq:implicit-nonconserved-energy-law}. The nonnegativity of
$\Phi$ then yields \eqref{eq:implicit-nonconserved-energy-decay}.
\end{proof}

\paragraph{Conserved variables}
For conserved variables, the network represents the flux
$\bm j_{\bm\theta}$ and is assumed to satisfy the prescribed
periodic or no-flux boundary conditions. The trial state is obtained
from the discrete continuity equation
\begin{equation*}
    \bm u^{k+1}_{\tau}(\bm\theta)
    =\bm u^k-\tau\nabla\cdot\bm j_{\bm\theta},
\end{equation*}
where the divergence acts componentwise. The corresponding implicit
variational problem is
\begin{equation}
    \bm\theta^{k+1}\in\arg\min_{\bm\theta\in\Theta_{\mathrm{ad}}^k}
    \left\{
       \mathcal E\bigl(\bm u^{k+1}_{\tau}(\bm\theta)\bigr)
       +\tau\Phi(\bm u^k;\bm j_{\bm\theta})
    \right\},
    \label{eq:implicit-conserved-minimization}
\end{equation}
and the accepted flux determines the new state through
\begin{equation}
    \bm u^{k+1}
    =\bm u^k-\tau\nabla\cdot\bm j_{\bm\theta^{k+1}}.
    \label{eq:implicit-conserved-update}
\end{equation}
Thus, the state and flux are linked by a conservative update, while
the new flux is selected through the discrete Onsager functional.

\begin{proposition}[Energy dissipation and mass conservation]
\label{prop:implicit-conserved-structure}
Let $\tau>0$, and suppose that the zero neural flux is admissible,
$\bm j_{\bm\theta_0}=\bm0$. Assume that
$\Phi(\bm u;\bm0)=0$ and $\Phi(\bm u;\bm j)\geq0$
for every admissible flux. 
For any admissible parameter $\bm\theta^{k+1}$ whose
neural flux satisfies the prescribed periodic or no-flux
boundary conditions, the update
\eqref{eq:implicit-conserved-update} conserves the mass
of each species:
\begin{equation}
    \int_\Omega u_i^{k+1}\,\mathrm{d}\bm x
    =
    \int_\Omega u_i^k\,\mathrm{d}\bm x,
    \qquad i=1,\ldots,s.
    \label{eq:implicit-conserved-mass}
\end{equation}

If, in addition, $\bm\theta^{k+1}$, not necessarily
a minimizer of \eqref{eq:implicit-conserved-minimization},
satisfies
\begin{equation}
    \mathcal J_k(\bm\theta^{k+1})
    \leq
    \mathcal J_k(\bm\theta_0)
    =
    \mathcal E(\bm u^k),
    \label{eq:implicit-conserved-descent}
\end{equation}
then the following discrete energy-dissipation
inequality holds:
\begin{equation}
    \mathcal E(\bm u^{k+1})
    +
    \tau\Phi(\bm u^k;\bm j_{\bm\theta^{k+1}})
    \leq
    \mathcal E(\bm u^k),
    \label{eq:implicit-conserved-energy-law}
\end{equation}
and consequently
\begin{equation}
    \mathcal E(\bm u^{k+1})
    \leq
    \mathcal E(\bm u^k).
    \label{eq:implicit-conserved-energy-decay}
\end{equation}
\end{proposition}

\begin{proof}
Integrating the conservative update
\eqref{eq:implicit-conserved-update} over $\Omega$
and applying the divergence theorem gives
\begin{align*}
    \int_\Omega (u_i^{k+1}-u_i^k)\,\mathrm{d}\bm x
    &=
    -\tau\int_\Omega
    \nabla\cdot\bm j_{i,\bm\theta^{k+1}}\,\mathrm{d}\bm x\\
    &=
    -\tau\int_{\partial\Omega}
    \bm j_{i,\bm\theta^{k+1}}\cdot\bm n\,\mathrm{d}S
    =0,
\end{align*}
where the boundary integral vanishes under the no-flux
condition or by cancellation under periodic boundary
conditions. This establishes
\eqref{eq:implicit-conserved-mass}, independently
of the optimization accuracy.

Furthermore, by the definition of $\mathcal J_k$
and the update \eqref{eq:implicit-conserved-update},
the descent condition
\eqref{eq:implicit-conserved-descent} yields
\eqref{eq:implicit-conserved-energy-law}.
The nonnegativity of $\Phi$ then implies
\eqref{eq:implicit-conserved-energy-decay}.
\end{proof}

\begin{remark}[Energy stability under inexact optimization]
\label{rem:implicit-inexact-optimization}
The global-minimizer assumption in the ideal implicit problems
\eqref{eq:implicit-nonconserved-minimization} and
\eqref{eq:implicit-conserved-minimization} is sufficient but not
necessary for energy dissipation. In either case, the zero process
variable satisfies
$\mathcal J_k(\bm\theta_0)=\mathcal E(\bm u^k)$.
After any finite number of optimization iterations, let
$\widetilde{\bm\theta}^{k+1}$ be the candidate produced by the
optimizer. The accepted parameter is defined by
\begin{equation}
    \bm\theta^{k+1}=
    \begin{cases}
      \widetilde{\bm\theta}^{k+1},
        &\widetilde{\bm\theta}^{k+1}\in\Theta_{\mathrm{ad}}^k
         \text{ and }
         \mathcal J_k(\widetilde{\bm\theta}^{k+1})
           \leq\mathcal J_k(\bm\theta_0),\\[2mm]
      \bm\theta_0,&\text{otherwise}.
    \end{cases}
    \label{eq:implicit-safeguarded-selection}
\end{equation}
This acceptance criterion enforces the descent conditions of both
propositions, irrespective of whether the nonlinear optimization
converges. Hence, energy dissipation is unconditional with respect to
$\tau$, although a rejected step may yield a zero update and poor
optimization may reduce accuracy or stall the evolution. When the
energy and objective are evaluated numerically, the same argument
applies to the consistently computed discrete energy (for example,
with a fixed quadrature); decay of the exact continuous energy
requires additional control of the integration error.
\end{remark}

Algorithm~\ref{alg:i-onts} gives the I-ONTS procedure in the unified
notation of Section~3.1. In practice, a gradient-based optimizer is
used to obtain a candidate for \eqref{eq:implicit-unified-objective};
the acceptance check is performed before the state is committed.

\begin{algorithm}[H]
\caption{Implicit Onsager neural time-stepping (I-ONTS)}
\label{alg:i-onts}
\hspace*{0.02in} {\bf Require:}
        Initial state $\bm u^0$; zero-process parameter $\bm\theta_0$; 
    Step size $\tau>0$; number of steps $N$; network $\bm z_{\bm\theta}$; 
    kinematic operator $\mathcal B$; Free energy $\mathcal E$ 
    and dissipation potential $\Phi$.
\begin{algorithmic}[1]
\State Compute an Onsager-consistent target $\bm z_{\mathrm{tar}}^0$ at $\bm u^0$.
\State Fit $\bm\theta^0$ via neural network training so that $\bm z_{\bm\theta^0}\approx\bm z_{\mathrm{tar}}^0$.
\For{$k=0,1,\ldots,N-1$}
    \State Define $\mathcal J_k(\bm\theta)$ by 
    \eqref{eq:implicit-unified-objective} and approximately minimize $\mathcal J_k$ starting from
           $\bm\theta^k$ to obtain $\widetilde{\bm\theta}^{k+1}$.
    \If{$\widetilde{\bm\theta}^{k+1}\in\Theta_{\mathrm{ad}}^k$
         \textbf{and} $\mathcal J_k(\widetilde{\bm\theta}^{k+1})
         \leq\mathcal E(\bm u^k)$}
         \State $\bm\theta^{k+1}\gets\widetilde{\bm\theta}^{k+1}$,
    \Else
         \State $\bm\theta^{k+1}\gets\bm\theta_0$.
    \EndIf
    \State $\bm u^{k+1}\gets\bm u^k+
           \tau\mathcal B\bm z_{\bm\theta^{k+1}}$,
\EndFor
\end{algorithmic}
\hspace*{0.02in} {\bf Output:} States $\{\bm u^k\}_{k=0}^{N}$ and accepted process parameters $\{\bm\theta^k\}_{k=0}^{N}$
\end{algorithm}

\subsection{Predictor--corrector Onsager neural time-stepping (PC-ONTS)}

PC-ONTS combines the explicit least-squares update of E-ONTS with
the variational correction of I-ONTS. The explicit stage constructs
an initial guess for the nonlinear minimization, while the corrected
state is always computed from the unchanged previous state $\bm u^k$.
For a given pair $(\bm u^k,\bm\theta^k)$, the explicit predictor
first forms the intermediate state
\begin{equation}
    \bm u_{\mathrm{pred}}^{k+1}
    :=\bm u^k+\tau\mathcal B\bm z_{\bm\theta^k}.
    \label{eq:pc-predictor-state}
\end{equation}
An Onsager-consistent target process variable is computed at
$\bm u_{\mathrm{pred}}^{k+1}$ and projected onto the linearized
network representation, exactly as in E-ONTS. This gives the
predictor parameters $\bm\theta_{\mathrm{LS}}^{k+1}$. The intermediate
state is not accepted as the final state; it is used only to
construct the predictor for the subsequent correction.

For multidimensional problems ($d\geq2$), we employ
an additional Rayleighian-based safeguard to improve
the reliability of the predictor initialization. 
At the fixed state $\bm u^k$,
we define
\begin{equation}
    \begin{aligned}
    \mathcal R_{\tau}^{k}(\bm u^k;\bm\theta)
    &:=
    \frac{
        \mathcal E\bigl(\bm u_\tau^{k+1}(\bm\theta)\bigr)
        -\mathcal E(\bm u^k)
    }{\tau}
    +\Phi(\bm u^k;\bm z_{\bm\theta})\\
    &=
    \frac{
        \mathcal J_k(\bm\theta)-\mathcal E(\bm u^k)
    }{\tau}.
    \end{aligned}
    \label{eq:pc-discrete-rayleighian}
\end{equation}
The least-squares predictor is used to initialize the implicit
corrector only if
\begin{equation}
    \mathcal R_{\tau}^{k}
       (\bm u^k;\bm\theta_{\mathrm{LS}}^{k+1})
    \leq
    \mathcal R_{\tau}^{k}(\bm u^k;\bm\theta^{k}).
    \label{eq:Rayleighian-based_safeguard}
\end{equation}
Otherwise, $\bm\theta^k$ is retained as the initial guess. 
For one-dimensional problems, the least-squares predictor is used
directly without this comparison. When state constraints are present,
only admissible initialization candidates are used; an inadmissible
predictor or old parameter is replaced by a feasible zero-process
initialization if necessary. The complete PC-ONTS procedure is summarized in
Algorithm~\ref{alg:pc-onts}.

The Rayleighian safeguard
\eqref{eq:Rayleighian-based_safeguard}
is used only to select the initialization for the
implicit corrector and does not replace the final
energy-descent acceptance criterion.
Consequently, under the assumptions of
Propositions~\ref{prop:implicit-nonconserved-energy}
and~\ref{prop:implicit-conserved-structure},
PC-ONTS inherits the discrete energy-dissipation
property of I-ONTS.
For conserved systems, the corrected update also preserves
mass under the prescribed periodic or no-flux
boundary conditions.

The three schemes offer complementary choices for updating neural
Onsager process variables. E-ONTS avoids nonlinear variational
minimization but does not generally guarantee energy decay;
I-ONTS retains the discrete Onsager structure; and PC-ONTS uses an
explicit predictor to initialize the safeguarded implicit step.
In the following section, we assess their accuracy, computational
cost, mass conservation, and energy-dissipation behavior on
representative gradient-flow PDEs, including increasingly stiff
regimes.

\begin{algorithm}[H]
\caption{Predictor--corrector Onsager neural time-stepping (PC-ONTS)}
\label{alg:pc-onts}
\hspace*{0.02in} {\bf Require:}
        Initial state $\bm u^0$; zero-process parameter
    $\bm\theta_0$; step size $\tau>0$; spatial dimension $d$;
    number of steps $N$; process network $\bm z_{\bm\theta}$;
    kinematic operator $\mathcal B$; free energy $\mathcal E$
    and dissipation potential $\Phi$.
\begin{algorithmic}[1]
\State Compute an Onsager-consistent target $\bm z_{\mathrm{tar}}^0$ at $\bm u^0$.
\State Fit $\bm\theta^0$ via neural network training so that $\bm z_{\bm\theta^0}\approx\bm z_{\mathrm{tar}}^0$.
\For{$k=0,1,\ldots,N-1$}
    \State $\bm u^{k+1}_{\mathrm{pred}}\gets\bm u^k+\tau\mathcal B\bm z_{\bm\theta^k}$, $\bm\mu_{\mathrm{pred}}^{k+1}\gets\delta\mathcal E(\bm u_{\mathrm{pred}}^{k+1})/\delta\bm u$
    \State Determine the Onsager target $\displaystyle
      \bm z_{\mathrm{tar}}^{k+1}\in
      \arg\min_{\bm z\in\mathcal Z_{\mathrm{ad}}}
      \left\{\langle\bm\mu^{k+1}_{\mathrm{pred}},\mathcal B\bm z\rangle
      +\Phi(\bm u_{\mathrm{pred}}^{k+1};\bm z)\right\}$

    \State Let $\bm J^k= D_{\bm\theta}\bm z_{\bm\theta^k}$ and determine the parameter velocity
    $\displaystyle
      \bm\eta^k\in\arg\min_{\bm\eta\in\mathbb R^p}
      \left\|\bm J^k\bm\eta
      -\frac{\bm z_{\mathrm{tar}}^{k+1}-\bm z_{\bm\theta^k}}{\tau}
      \right\|_{L^2(\Omega)}^2$
        \State $\bm\theta_{\mathrm{LS}}^{k+1}\gets\bm\theta^k+\tau\bm\eta^k$
        \If{$d\geq2$}
            \If{$\bm\theta_{\mathrm{LS}}^{k+1}$ is admissible
                 \textbf{and}
                 $\mathcal R_{\tau}^{k}
                     (\bm u^k;\bm\theta_{\mathrm{LS}}^{k+1})
                    \leq\mathcal R_{\tau}^{k}(\bm u^k;\bm\theta^k)$}
                \State $\bm\theta_{\mathrm{init}}
                       \gets\bm\theta_{\mathrm{LS}}^{k+1}$
            \Else
                \State $\bm\theta_{\mathrm{init}}\gets\bm\theta^k$
            \EndIf
        \Else
            \State $\bm\theta_{\mathrm{init}}
                  \gets\bm\theta_{\mathrm{LS}}^{k+1}$
        \EndIf
    \State If needed, replace an inadmissible
           $\bm\theta_{\mathrm{init}}$ by $\bm\theta_0$
    \State Approximately minimize $\mathcal J_k$ \eqref{eq:implicit-unified-objective} starting from
           $\bm\theta_{\mathrm{init}}$ to obtain
           $\widetilde{\bm\theta}^{k+1}$
    \If{$\widetilde{\bm\theta}^{k+1}\in\Theta_{\mathrm{ad}}^k$
         \textbf{and}
         $\mathcal J_k(\widetilde{\bm\theta}^{k+1})
         \leq\mathcal E(\bm u^k)$}
        \State $\bm\theta^{k+1}
               \gets\widetilde{\bm\theta}^{k+1}$
    \Else
        \State $\bm\theta^{k+1}\gets\bm\theta_0$
    \EndIf
    \State $\bm u^{k+1}\gets
           \bm u^k+\tau\mathcal B\bm z_{\bm\theta^{k+1}}$
\EndFor
\end{algorithmic}
\hspace*{0.02in} {\bf Output:} States $\{\bm u^k\}_{k=0}^{N}$ and accepted process
        parameters $\{\bm\theta^k\}_{k=0}^{N}$.
\end{algorithm}

\section{Numerical Experiments}
\subsection{Experimental Setup}

All experiments are implemented in JAX, using double precision unless
otherwise stated. Hardware, software, and other reproducibility details
are provided in \ref{app:common_setup}. 
To measure the overall accuracy over the complete space-time domain, we use
the relative $L^2$ error
\begin{equation}
\label{eq:spacetime_error}
    \mathrm{Err}_{st}
    =
    \frac{
        \left\|u_{\mathrm{num}}-u_{\mathrm{ref}}\right\|_
        {L^2([0,T]\times\Omega)}
    }{
        \left\|u_{\mathrm{ref}}\right\|_
        {L^2([0,T]\times\Omega)}
    },
\end{equation}
where $u_{\mathrm{ref}}$ denotes either the exact solution, when available,
or a sufficiently accurate reference solution. For Monte Carlo
discretizations in high-dimensional examples, the corresponding spatial
$L^2$ integrals are approximated using independent test samples, as
specified in the individual experiments. 
For conservative problems, we monitor the deviation of the numerical total
mass from its initial value. We denote the numerical mass by
\begin{equation}
\label{eq:numerical_mass}
    M_h(t)
    =
    \int_{\Omega} u_h(x,t)\,dx,
\end{equation}
where $u_h$ denotes the numerical solution. The absolute mass error is then
defined as
\begin{equation}
\label{eq:mass_error}
    \mathrm{Err}_{M}(t)
    =
    \left|M_h(t)-M_h(0)\right|.
\end{equation}
For weighted density-ratio formulations, the physical mass is used.
The corresponding definition and the computational timing conventions
are also given in \ref{app:common_setup}.

\subsection{Accuracy and Temporal Convergence: 1D Heat Equation}
\label{experiment_1dheat}
We first consider a one-dimensional heat equation
with periodic boundary conditions as a basic benchmark
for assessing the accuracy and temporal convergence
of the proposed methods:
\begin{equation}
\label{1d_heat_eq}
    \begin{cases}
    \partial_t u(x,t)=\partial_{xx}u(x,t),
        & x\in(-l,l),\quad t>0,\\
    u(-l,t)=u(l,t),\quad
    \partial_xu(-l,t)=\partial_xu(l,t),
        & t\geq0,\\
    u(x,0)=u_0(x),
        & x\in[-l,l].
    \end{cases}
\end{equation}
We set $l=\pi$ and $u_0(x)=2+\sin x$.
The corresponding exact solution is
\begin{equation*}
    u_{\mathrm{ref}}(x,t)
    =
    2+e^{-t}\sin x.
\end{equation*}
The availability of an exact solution makes this example particularly suitable 
for quantitatively examining the temporal accuracy of EDNN, E-ONTS, I-ONTS and PC-ONTS. 
The computations use $N_x=1024$ collocation points and $T=1$;
the network architecture, periodic embedding, and training settings
are detailed in \ref{app:one_d_settings}. 
Table~\ref{tab:heat_accuracy} reports the space-time relative 
errors for different time-step sizes. The amortized evolution cost 
per time step is also reported for each method.

As shown in Table~\ref{tab:heat_accuracy}, for this 
simple one-dimensional heat equation, the proposed E-ONTS, I-ONTS, 
and PC-ONTS methods do not exhibit a pronounced advantage over the 
classical EDNN method in terms of accuracy. All four methods achieve 
comparable space-time relative $L^2$ errors, which is expected 
since the present problem is simple and only mildly stiff. The main purpose of 
this experiment is therefore to verify the effectiveness and temporal 
convergence behavior of the proposed methods in a basic setting.
The errors of all four methods decrease at approximately first-order 
rates as the time step is refined. This behavior is consistent with the temporal 
discretization adopted here: the evolution equations are advanced 
using a first-order Euler scheme, while the spatial dependence of 
the solution is represented by neural networks. Higher-order time-stepping 
schemes \cite{du2021evolutional} could in principle be incorporated into the proposed framework 
to achieve higher temporal accuracy. However, the development of high-order 
temporal discretizations is not the focus of the present work; instead, 
we concentrate on the stability properties induced by different formulations 
of the neural evolution equations.

\begin{table}[H]
\centering
\caption{Space-time relative $L^2$ errors and mean evolution cost per time step for the one-dimensional heat equation.}
\label{tab:heat_accuracy}
\resizebox{\textwidth}{!}{
\begin{tabular}{c|cc|cc|cc|cc}
\toprule
& \multicolumn{2}{c|}{EDNN}
& \multicolumn{2}{c|}{E-ONTS}
& \multicolumn{2}{c|}{I-ONTS}
& \multicolumn{2}{c}{PC-ONTS}\\
\(\Delta t\)
& Error & Time/step (s)
& Error & Time/step (s)
& Error & Time/step (s)
& Error & Time/step (s)\\
\midrule

$4.0\times10^{-2}$
& $1.8881\times10^{-3}$ & $0.0916$
& $1.7645\times10^{-3}$ & $0.1221$
& $2.1393\times10^{-3}$ & $0.0589$
& $1.9236\times10^{-3}$ & $0.1613$\\

$2.0\times10^{-2}$
& $9.1084\times10^{-4}$ & $0.0917$
& $9.1722\times10^{-4}$ & $0.1208$
& $1.1155\times10^{-3}$ & $0.0444$
& $1.0684\times10^{-3}$ & $0.1519$\\

$1.0\times10^{-2}$
& $4.2661\times10^{-4}$ & $0.0897$
& $4.6473\times10^{-4}$ & $0.1251$
& $5.7352\times10^{-4}$ & $0.0520$
& $5.4895\times10^{-4}$ & $0.1489$\\

$5.0\times10^{-3}$
& $1.8676\times10^{-4}$ & $0.1080$
& $2.3371\times10^{-4}$ & $0.1184$
& $2.7914\times10^{-4}$ & $0.0818$
& $2.7084\times10^{-4}$ & $0.1630$\\

$2.5\times10^{-3}$
& $7.1191\times10^{-5}$ & $0.0904$
& $1.1756\times10^{-4}$ & $0.1170$
& $1.0853\times10^{-4}$ & $0.0555$
& $1.1303\times10^{-4}$ & $0.1459$\\

\midrule
\multicolumn{1}{c|}{Average time/step (s)}
& \multicolumn{2}{c|}{$0.0943$}
& \multicolumn{2}{c|}{$0.1207$}
& \multicolumn{2}{c|}{$0.0585$}
& \multicolumn{2}{c}{$0.1542$}\\
\bottomrule
\end{tabular}
}
\end{table}

\begin{figure}[H]
    \centering
    \begin{subfigure}{0.48\textwidth}
        \centering
        \includegraphics[width=\textwidth]{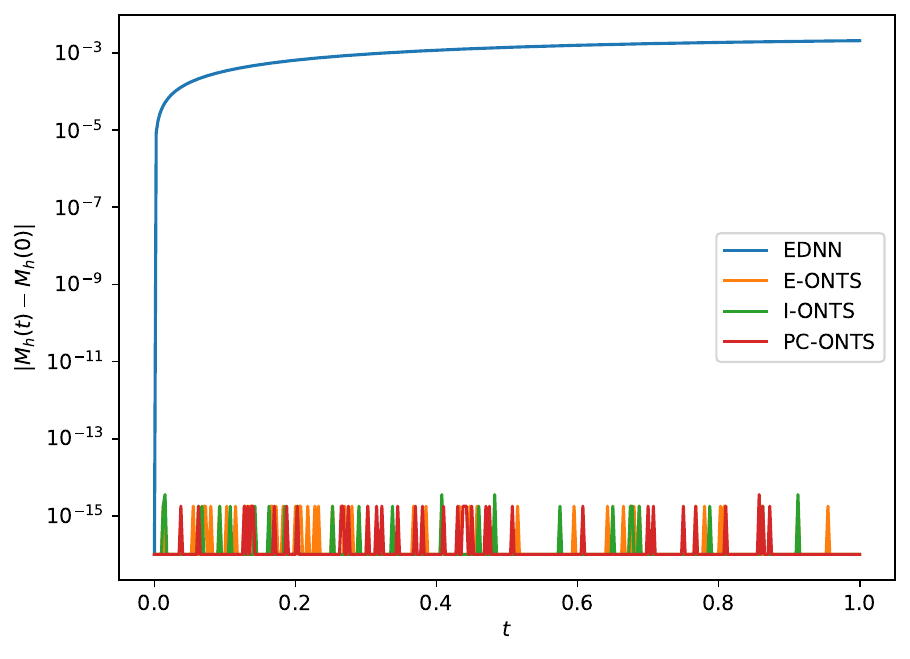}
        \caption{Absolute mass error.}
        \label{fig:heat_mass_error}
    \end{subfigure}
    \hfill
    \begin{subfigure}{0.48\textwidth}
        \centering
        \includegraphics[width=\textwidth]{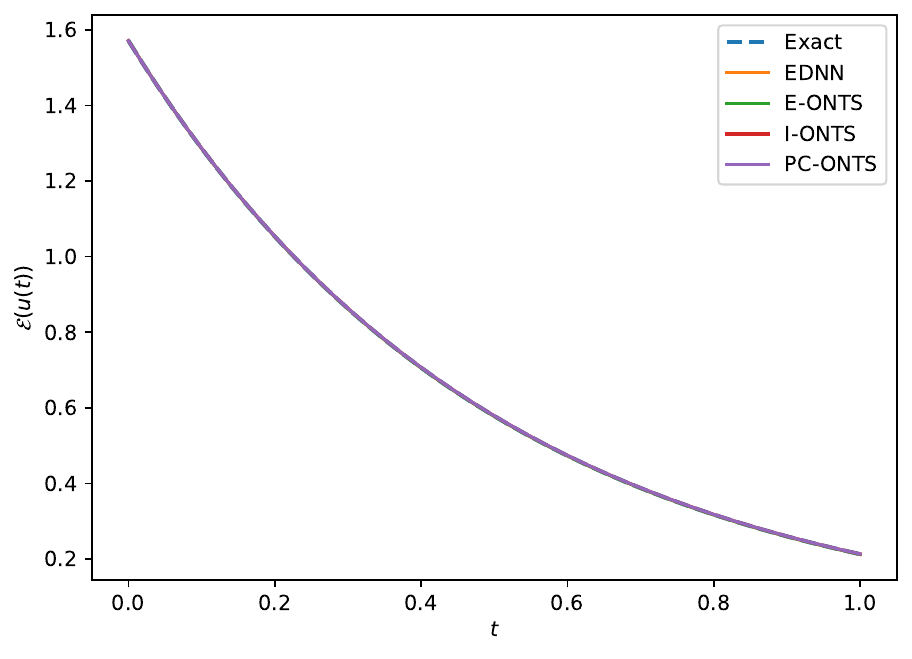}
        \caption{Free-energy evolution.}
        \label{fig:heat_energy_error}
    \end{subfigure}

    \caption{
    Mass conservation and free-energy dissipation for
    the one-dimensional heat equation with
    $\Delta t=2.5\times10^{-3}$:
    (a) absolute mass error;
    (b) free-energy evolution.
    }
    \label{fig:heat_mass_energy}
\end{figure}

As shown in Table~\ref{tab:heat_accuracy}, PC-ONTS incurs
the cost of both least-squares prediction and implicit correction,
whereas the other methods involve only their respective evolution
steps. A detailed discussion of per-step costs is provided in
\ref{app:one_d_settings}. We further examine the structural properties of the numerical 
solutions obtained with $\Delta t=2.5\times10^{-3}$. 
Figure~\ref{fig:heat_mass_energy} compares the absolute 
mass errors (left panel) and free-energy evolution (right panel) 
of the four methods. Consistent with the mass-conservation properties 
established in Propositions~\ref{prop:explicit-mass-conservation} 
and~\ref{prop:implicit-conserved-structure}, the mass errors of 
E-ONTS, I-ONTS, and PC-ONTS remain around $10^{-15}$, close to 
numerical round-off, whereas the mass error of EDNN grows to 
approximately $10^{-3}$. This demonstrates the advantage of the 
Onsager-based formulations in preserving mass. For this relatively 
simple, linear one-dimensional heat equation, all four methods also 
closely reproduce the exact free-energy decay, yielding nearly 
indistinguishable energy curves. Overall, the proposed methods 
achieve accuracy comparable to EDNN for this problem while preserving 
mass much more accurately. We next consider a one-dimensional 
Allen--Cahn equation to investigate the potential computational 
advantages of the predictor--corrector strategy in nonlinear gradient-flow problems.

\subsection{Efficiency of Predictor--Corrector Onsager Neural Time-Stepping}

We examine the efficiency of the proposed predictor--corrector
strategy using the one-dimensional Allen--Cahn equation:
\begin{equation}
\label{eq:allen_cahn}
    \begin{cases}
    \partial_t u(x,t)
    =
    \gamma_1 \partial_{xx}u(x,t)
    +
    \gamma_2\left(u(x,t)-u^3(x,t)\right),
    & x\in(-1,1),\quad t>0,\\
    u(-1,t)=u(1,t),\,
    \partial_xu(-1,t)=\partial_xu(1,t),
    & t\geq0,\\
    u(x,0)=u_0(x),
    & x\in[-1,1],
    \end{cases}
\end{equation}
with
\begin{equation*}
    u_0(x)
    =
    0.45\cos(\pi x)
    -0.25\sin(2\pi x)
    +0.18\cos(3\pi x)
    +0.10\sin(5\pi x),
\end{equation*}
and $\gamma_1=10^{-4},\,\gamma_2=5.$

The Allen--Cahn equation \eqref{eq:allen_cahn} is a nonconserved
gradient flow associated with the free energy
\begin{equation*}
    \mathcal E(u)
    =
    \int_{-1}^{1}
    \left[
        \frac{\gamma_1}{2}|\partial_xu|^2
        +
        \frac{\gamma_2}{4}(u^2-1)^2
    \right]\,dx.
\end{equation*}
Indeed,
\begin{equation*}
    \partial_tu
    =
    -\frac{\delta\mathcal E}{\delta u}
    =
    \gamma_1\partial_{xx}u
    +
    \gamma_2(u-u^3).
\end{equation*}
Accordingly, the Onsager-based methods represent the temporal rate
$v_{\bm{\theta}}\approx\partial_t u$.
We use $\Delta t=10^{-3}$ and $T=1$.
All methods employ periodic neural representations and the common
collocation and training settings specified in \ref{app:one_d_settings}.
Since a closed-form transient solution is unavailable,
a high-accuracy numerical solution is used as the reference.

Table~\ref{tab:allen_cahn_accuracy} reports the space-time relative $L^2$ errors of different methods. All three Onsager-based methods achieve substantially smaller errors than EDNN, with PC-ONTS attaining the highest accuracy using only 20 optimization iterations per time step, compared with 400 iterations for I-ONTS. To further assess the computational efficiency of the predictor--corrector strategy, we compare I-ONTS and PC-ONTS under different optimization budgets in Table~\ref{tab:allen_cahn_cost}.

\begin{table}[H]
\centering
\caption{Space-time relative $L^2$ errors for the 1D Allen--Cahn equation.}
\label{tab:allen_cahn_accuracy}
\begin{tabular}{c cccc}
\toprule
$\Delta t$ & EDNN & E-ONTS & I-ONTS & PC-ONTS \\
\midrule
$1.0\times10^{-3}$
& $3.8586\times10^{-2}$
& $1.1399\times10^{-3}$
& $1.1258\times10^{-3}$
& $\mathbf{6.4674\times10^{-4}}$\\
\bottomrule
\end{tabular}
\end{table}

\begin{table}[H]
\centering
\caption{Space-time relative $L^2$ errors and mean wall-clock
time per time step for I-ONTS and PC-ONTS on the 1D
Allen--Cahn equation. Here, $N_{\rm opt}$ denotes the 
number of neural-network optimization iterations for the implicit correction
at each time step.}
\label{tab:allen_cahn_cost}
\begin{tabular}{c|c|cc}
\toprule
Method & $N_{\rm opt}$ & Relative $L^2$ error & Time/step (s)\\
\midrule
I-ONTS  & $20$  & $1.3178\times10^{-2}$ & $0.0530$\\
I-ONTS  & $200$ & $1.7846\times10^{-3}$ & $0.1514$\\
I-ONTS  & $400$ & $1.1258\times10^{-3}$ & $0.2518$\\
PC-ONTS & $20$  & $\mathbf{6.4674\times10^{-4}}$ & $0.1669$\\
\bottomrule
\end{tabular}
\end{table}

\begin{figure}[H]
    \centering
    \begin{subfigure}{0.48\textwidth}
        \centering
        \includegraphics[width=\textwidth]
        {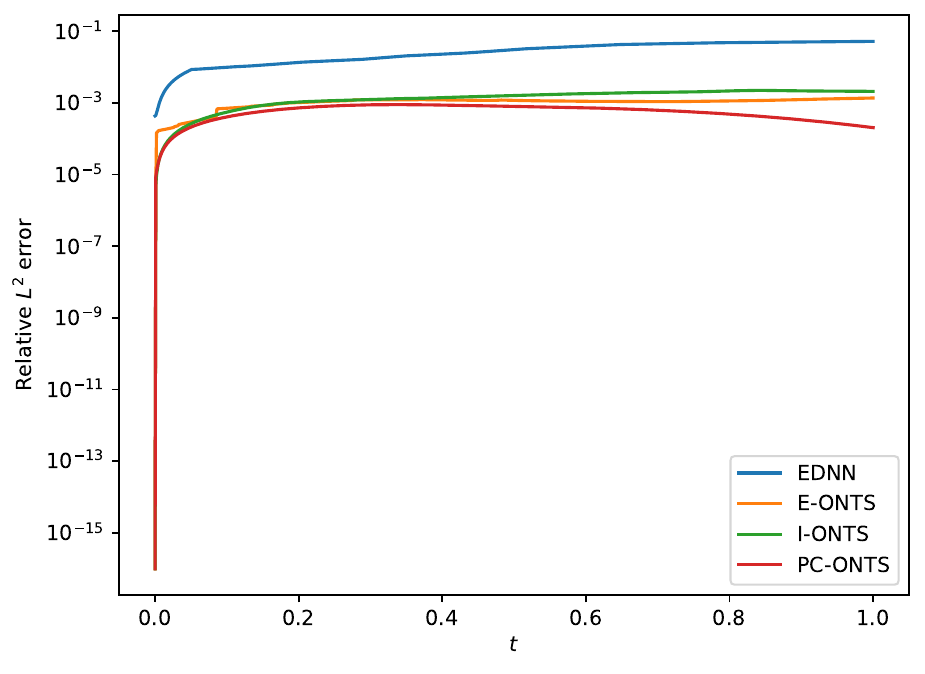}
        \caption{Relative $L^2$ error of the solution.}
        \label{fig:allen_cahn_error}
    \end{subfigure}
    \hfill
    \begin{subfigure}{0.48\textwidth}
        \centering
        \includegraphics[width=\textwidth]
        {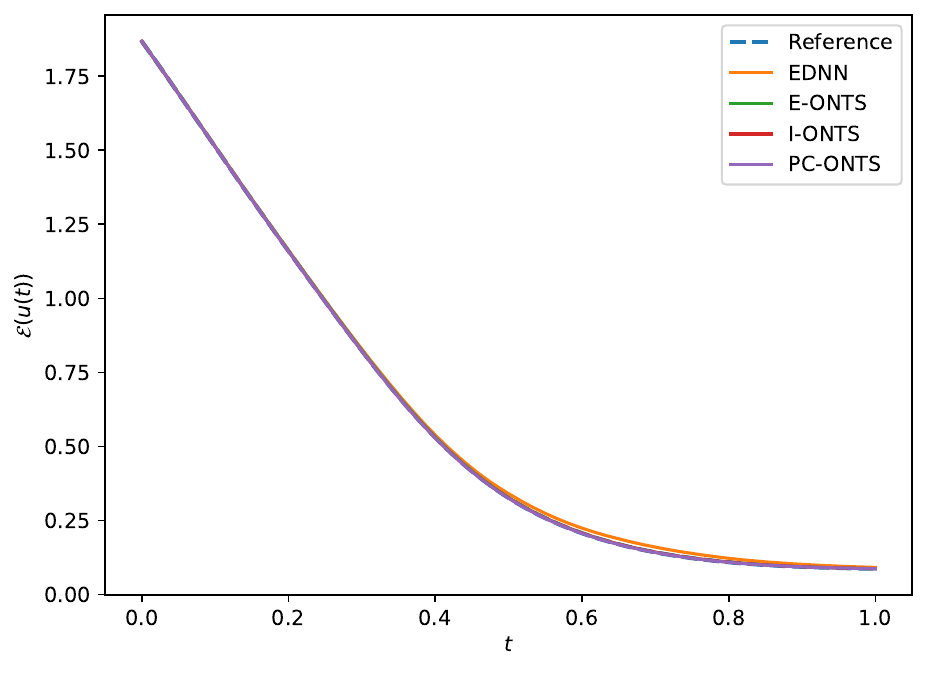}
        \caption{Free-energy evolution.}
        \label{fig:allen_cahn_energy}
    \end{subfigure}

    \caption{
    Accuracy and energy evolution for the one-dimensional
    Allen--Cahn equation with $\Delta t=10^{-3}$.
    }
    \label{fig:allen_cahn_accuracy_energy}
\end{figure}

The results show that PC-ONTS achieves an error of $6.4674\times10^{-4}$ with only 20 optimization iterations per time step, outperforming I-ONTS with 400 iterations, which yields an error of $1.1258\times10^{-3}$. Moreover, the mean computational cost per time step is reduced from $0.2518$ to $0.1669$ seconds. Although the least-squares predictor introduces additional computational overhead compared with I-ONTS using the same optimization budget, it substantially reduces the optimization effort required for the implicit correction. These results demonstrate that the predictor--corrector strategy can improve the accuracy--cost trade-off when the explicit prediction provides a reliable initialization. However, in higher-dimensional and more challenging problems, the predictor may be rejected more frequently, potentially reducing the efficiency gains of PC-ONTS. Nevertheless, this safeguard is essential for preventing unreliable predictions from being used in the implicit correction, thereby maintaining the reliability of the overall scheme.
Figure~\ref{fig:allen_cahn_accuracy_energy} illustrates the relative $L^2$ error and free-energy evolution. The error curves further demonstrate the improved accuracy of the proposed Onsager-based methods compared with EDNN, while the energy curves closely follow the reference evolution, confirming that all three schemes capture the underlying energy-dissipation behavior.
To further assess the stability and robustness of the proposed schemes under more demanding conditions, we next consider two-dimensional problems that pose greater numerical challenges.

\subsection{Stability and Robustness of Onsager Neural Time-Stepping}

The preceding experiments mainly focus on solution accuracy and structure
preservation. We now turn to the stability and robustness of the proposed schemes 
through two complementary tests. 
First, we consider a two-dimensional heat equation with varying diffusion coefficients, keeping the time step fixed while progressively increasing the diffusion coefficient. This
provides a controlled setting for examining the behavior of the methods as
the problem becomes increasingly stiff. We then consider a two-dimensional
Cahn--Hilliard equation, which involves nonlinear fourth-order dynamics and
places more stringent requirements on numerical stability. Together, these
two experiments assess the robustness of the proposed formulations under
both systematically increased stiffness and more challenging nonlinear
dynamics.

\subsubsection{Stiff Heat Equation in Two Dimensions}

We consider the two-dimensional periodic heat equation
\begin{equation}
\label{eq:heat2d_stability}
\begin{cases}
\partial_t u(x,y,t)
=
\kappa\left(
\partial_{xx}u(x,y,t)
+
\partial_{yy}u(x,y,t)
\right),
&
(x,y)\in(-\pi,\pi)^2,\quad t>0,
\\
u(x,y,0)
=
2+\sin(x)\sin(y)
+
a\sin(Kx)\sin(Ky),
&
(x,y)\in[-\pi,\pi]^2,
\end{cases}
\end{equation}
subject to periodic boundary conditions in both spatial directions.
Throughout this experiment, we fix
\begin{equation*}
    a=0.2,\qquad K=2,\qquad \Delta t=5\times10^{-2}, \qquad T=1,
\end{equation*}
and vary the diffusion coefficient $\kappa$.
The exact solution of \eqref{eq:heat2d_stability} is
\begin{equation*}
    u_{\mathrm{ref}}(x,y,t)
    =2+e^{-2\kappa t}\sin(x)\sin(y)+
    a e^{-2\kappa K^2t}\sin(Kx)\sin(Ky).
\end{equation*}
Thus, the high-frequency component decays on the time scale
$(2\kappa K^2)^{-1}$. 
Increasing $\kappa$ shortens this time scale and makes the 
problem increasingly stiff relative to the fixed time step, 
thereby imposing a more restrictive stability condition on 
explicit time discretizations. 
For every value of $\kappa$, the collocation grid, network
architectures, and optimization settings are kept fixed
(see \ref{app:two_d_settings}). 
The initial state is independent of $\kappa$, while the 
neural flux is optimized separately for each diffusion coefficient.

To understand the stability behavior of the four methods,
we first consider the classical forward-Euler discretization 
of the heat equation, which provides a useful reference for
the explicit temporal treatments employed in EDNN and
E-ONTS. Applying the forward-Euler method to a Fourier
mode with wave number $(K,K)$ gives the amplification
factor
\begin{equation}
\label{eq:heat2d_explicit_amplification}
G_{\rm E}
=
1-2\kappa K^2\Delta t.
\end{equation}
The corresponding stability condition is
\begin{equation}
\label{eq:heat2d_explicit_condition}
|G_{\rm E}|\leq1
\quad\Longleftrightarrow\quad
\kappa K^2\Delta t\leq1.
\end{equation}
We therefore introduce the dimensionless parameter
\begin{equation}
\label{eq:heat2d_stiffness_number}
\eta=\kappa K^2\Delta t
\end{equation}
to characterize the stiffness relative to the prescribed
time step. For $K=2$ and $\Delta t=0.05$, we have
$\eta=0.2\kappa$. Although this classical stability 
condition does not directly characterize the stability 
of the neural-network-based schemes, it provides a 
useful reference for understanding the increasing 
stability challenges faced by explicit temporal 
treatments as $\kappa$ grows. 

In contrast, I-ONTS employs an implicit Onsager
variational update, while PC-ONTS combines an explicit
predictor with the same implicit variational correction.
Both methods are based on the unified objective
\eqref{eq:implicit-unified-objective}. For the present
heat equation, the free energy and dissipation potential
are given by
\begin{equation*}
    \mathcal{E}(u)
    =
    \frac{1}{2}\int_{\Omega}(u-2)^2\,\mathrm{d}\bm{x},
    \qquad
    \Phi(u;\bm{j})
    =
    \frac{1}{2\kappa}
    \int_{\Omega}|\bm{j}|^2\,\mathrm{d}\bm{x}.    
\end{equation*}
Using the neural flux representation, the corresponding
implicit objective becomes
\begin{equation*}
    \mathcal{J}_k(\bm{\theta})
    =
    \frac{1}{2}
    \int_{\Omega}
    \left(
    u_{\Delta t}^{k+1}(\bm{\theta})-2
    \right)^2\,\mathrm{d}\bm{x}
    +
    \frac{\Delta t}{2\kappa}
    \int_{\Omega}
    |\bm{j}_{\bm{\theta}}|^2\,\mathrm{d}\bm{x},
\end{equation*}
where
\begin{equation*}
    u_{\Delta t}^{k+1}(\bm{\theta})=u^k-\Delta t\,\nabla\cdot\bm{j}_{\bm{\theta}}.
\end{equation*}
For comparison, we consider the corresponding
variational minimization over the full admissible
flux space. The first-order optimality condition yields
$\bm{j}=-\kappa\nabla u^{k+1}$, which recovers 
the standard implicit-Euler discretization of the heat equation.
For the Fourier mode with wave number $(K,K)$,
the corresponding amplification factor is
\begin{equation}
\label{eq:heat2d_implicit_amplification}
    G_{\rm I}
    =
    \frac{1}{1+2\kappa K^2\Delta t},
\end{equation}
which satisfies $0<G_{\rm I}\leq1$ for any
$\kappa\Delta t>0$. 
This comparison highlights the stability advantage of
the underlying implicit variational discretization,
which is shared by I-ONTS and the corrector step of
PC-ONTS. Although neural-network approximation and 
inexact optimization may affect the actual numerical 
performance, the implicit variational structure provides 
a basis for unconditional discrete energy stability, 
subject to the suitable descent criterion established above. 
Then, We examine whether these theoretical stability 
characteristics are reflected in the numerical performance of the four methods. 
Table~\ref{tab:heat2d_stability} reports the
space-time relative $L^2$ errors of the four methods
as $\kappa$ increases with $\Delta t$ fixed.

Table~\ref{tab:heat2d_stability} reveals increasingly pronounced 
differences in the numerical performance of the four methods as $\kappa$
increases. For $\kappa=1$ and $2$, all methods remain numerically stable,
although E-ONTS already exhibits a substantial increase in error at $\kappa=2$
and breaks down for $\kappa\geq4$. EDNN remains stable over a wider range,
but its error increases from $8.22\times10^{-3}$ at $\kappa=1$ to $6.93\times10^{-1}$
at $\kappa=32$, before breaking down at $\kappa=64$. In contrast, I-ONTS and PC-ONTS,
both employing an implicit Onsager variational correction, remain stable throughout the
tested range. While the error of I-ONTS reaches approximately $10^{-1}$ at $\kappa=32$
and $64$, PC-ONTS maintains errors around $10^{-2}$ or below and achieves the highest
accuracy for all tested values of $\kappa$. The improved performance of PC-ONTS may be
attributed in part to the explicit least-squares predictor and Rayleighian-based
safeguard \eqref{eq:Rayleighian-based_safeguard}, which together provide a more 
suitable initialization for the implicit correction.
These results are qualitatively consistent with the preceding modal analysis,
while also reflecting the effects of neural-network approximation and
optimization on numerical stability.

\begin{table}[H]
\centering
\caption{
Space-time relative $L^2$ errors for the two-dimensional heat equation
with $a=0.2$, $K=2$, and $\Delta t=5\times10^{-2}$.
The stiffness indicator is
$\eta=\kappa K^2\Delta t=0.2\kappa$.
``NaN/Inf'' indicates numerical breakdown during time evolution.
}
\label{tab:heat2d_stability}
\begin{tabular}{cccccc}
\toprule
$\kappa$
& $\eta$
& EDNN
& E-ONTS
& I-ONTS
& PC-ONTS
\\
\midrule

$1$
& $0.2$
& $8.2213\times10^{-3}$
& $4.6558\times10^{-3}$
& $4.4415\times10^{-3}$
& $\mathbf{3.9987\times10^{-3}}$
\\

$2$
& $0.4$
& $2.0084\times10^{-2}$
& $7.0979\times10^{-2}$
& $7.1320\times10^{-3}$
& $\mathbf{6.0720\times10^{-3}}$
\\

$4$
& $0.8$
& $5.0192\times10^{-2}$
& \textit{NaN/Inf}
& $1.4320\times10^{-2}$
& $\mathbf{8.1312\times10^{-3}}$
\\

$8$
& $1.6$
& $1.5277\times10^{-1}$
& \textit{NaN/Inf}
& $3.2037\times10^{-2}$
& $\mathbf{1.0383\times10^{-2}}$
\\

$16$
& $3.2$
& $4.4688\times10^{-1}$
& \textit{NaN/Inf}
& $4.8905\times10^{-2}$
& $\mathbf{1.2171\times10^{-2}}$
\\

$32$
& $6.4$
& $6.9278\times10^{-1}$
& \textit{NaN/Inf}
& $9.9727\times10^{-2}$
& $\mathbf{1.1678\times10^{-2}}$
\\

$64$
& $12.8$
& \textit{NaN/Inf}
& \textit{NaN/Inf}
& $ 1.0101\times10^{-1}$
& $\mathbf{7.4992\times10^{-3}}$
\\

\bottomrule
\end{tabular}
\end{table}

\subsubsection{Stiff Cahn--Hilliard Equation in Two Dimensions}

We further examine the robustness of the temporal discretizations using
the two-dimensional Cahn--Hilliard equation
\begin{equation}
\label{eq:ch2d_stability}
    \partial_t\phi
    =
    \lambda\Delta\mu,
    \qquad
    \mu
    =
    -\epsilon^2\Delta\phi
    +
    \phi^3-\phi,
    \qquad
    (x,y)\in[-1,1)^2,
\end{equation}
subject to periodic boundary conditions. We set $\lambda=1, \, \epsilon=0.05$ and
the initial condition consists of two diffuse bubbles. Let
\begin{equation*}
    d_{\pm}(x,y)
    =
    \sqrt{(x\mp0.28)^2+y^2}.
\end{equation*}
Taking the bubble radius as $r=0.4$, the initial phase field is defined by
\begin{equation*}
    \phi_0(x,y)
    =
    \max\left\{
    \tanh\left(
        \frac{r-d_{+}(x,y)}{2\epsilon}
    \right),
    \tanh\left(
        \frac{r-d_{-}(x,y)}{2\epsilon}
    \right)
    \right\}.
\end{equation*}
The associated free energy of \eqref{eq:ch2d_stability} is
\begin{equation}
\label{eq:ch2d_energy}
    \mathcal E(\phi)
    =
    \int_{\Omega}
    \left[
        \frac{\epsilon^2}{2}|\nabla\phi|^2
        +
        \frac{1}{4}(\phi^2-1)^2
    \right]\,\mathrm{d}\bm{x},
\end{equation}
and the conservative flux is
\begin{equation}
\label{eq:ch2d_flux}
    \bm{j}
    =-\lambda\nabla\mu,\qquad
    \partial_t\phi+\nabla\cdot\bm{j}=0.
\end{equation}

\begin{figure}[H]
    \centering
    \includegraphics[width=0.98\textwidth]
    {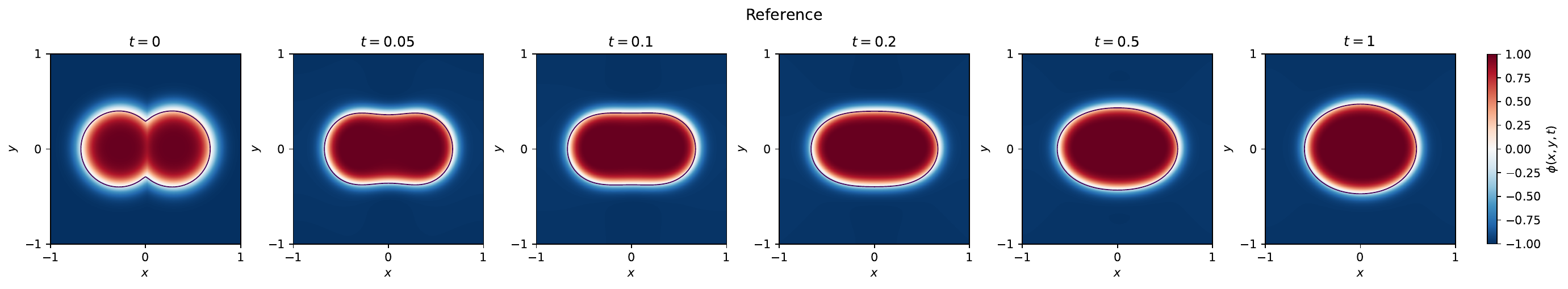}
    \caption{
    Reference phase-field evolution for the two-dimensional
    Cahn--Hilliard two-bubble problem.
    The black curves indicate the zero level set $\phi=0$,
    highlighting the evolution of the diffuse interfaces.
    }
    \label{fig:ch2d_reference_phase}
\end{figure}

The neural solvers use a $64\times64$ periodic collocation grid,
$\Delta t=10^{-3}$, and $T=1$. The fourth-order term introduces 
strong high-frequency stiffness, while the pointwise maximum 
used to construct the initial condition introduces a loss of 
smoothness along the symmetry line where the two bubble 
profiles coincide. To initialize
I-ONTS and PC-ONTS, the flux networks are fitted to the
first-step flux $\bm j^1$ at $t=\Delta t$, extracted from 
the Fourier pseudospectral reference trajectory. This 
reference-informed initialization provides a more 
suitable initial guess for the stiff dynamics and 
is used only for the first time step; all subsequent 
updates are performed using the respective neural 
time-stepping schemes. Its motivation, the network
and optimization settings, and the reference-solution protocol
are described in \ref{app:two_d_settings}.
The reference evolution, computed on a $256\times256$ grid with
$\Delta t_{\rm ref}=10^{-4}$, is illustrated in
Figure~\ref{fig:ch2d_reference_phase}.

Table~\ref{tab:ch2d_stability} reveals pronounced differences in the 
numerical performance of the four methods. Although EDNN produces
finite numerical solutions, its space-time relative $L^2$ error
exceeds $1$, while E-ONTS breaks down during the evolution. In
contrast, I-ONTS and PC-ONTS remain stable, achieving comparable
errors of $7.2627\times10^{-4}$ and $7.5129\times10^{-4}$, respectively.
Due to the stiffness associated with the fourth-order dynamics, the
explicit least-squares predictor is often unreliable in this experiment. 
Consequently, the Rayleighian-based safeguard \eqref{eq:Rayleighian-based_safeguard} 
rejects unreliable predictions and falls back to the previous parameter vector to 
initialize the implicit correction. 
This behavior is consistent with the similar accuracy achieved by I-ONTS and PC-ONTS.
Supplementary phase-field snapshots and pointwise error maps for
I-ONTS and PC-ONTS are shown in Figures~\ref{fig:ch2d_ionts_phase_error}
and~\ref{fig:ch2d_pconts_phase_error} of
\ref{app:two_d_settings}. Both methods accurately capture the reference
phase-field evolution, with small errors primarily concentrated near the moving interfaces.

\begin{table}[H]
\centering
\caption{
Space-time relative $L^2$ errors for the two-dimensional
Cahn--Hilliard problem.
``NaN/Inf'' indicates numerical breakdown during time evolution.
}
\label{tab:ch2d_stability}
\begin{tabular}{c|cccc}
\toprule
Method
& EDNN
& E-ONTS
& I-ONTS
& PC-ONTS
\\
\midrule
Space-time relative $L^2$ error
& $1.2346$
& \textit{NaN/Inf}
& $\mathbf{7.2627\times10^{-4}}$
& $7.5129\times10^{-4}$
\\
\bottomrule
\end{tabular}
\end{table}

Together with the preceding heat-equation experiment, these results 
demonstrate the robustness of the implicit Onsager formulations under 
challenging dynamics. In particular, PC-ONTS retains the stability 
of the implicit correction while exploiting the explicit predictor 
whenever it provides a reliable initialization. We next investigate 
the performance of the proposed schemes in high-dimensional and 
coupled systems, which introduce additional challenges in terms 
of computational complexity and numerical robustness.

\subsection{High-Dimensional and Coupled Systems}

We finally investigate the performance of the proposed schemes for 
high-dimensional and coupled PDE systems. As the spatial dimension 
increases, conventional grid-based discretizations become increasingly 
expensive. Neural-network-based spatial representations combined with 
Monte Carlo sampling provide an alternative that avoids the exponential 
growth associated with tensor-product grids. 
We consider high-dimensional Allen--Cahn and Fokker--Planck 
equations as representative examples of nonconserved and conserved 
dynamics, respectively, followed by a coupled multi-species Poisson--Nernst--Planck 
system involving electrostatic interactions. These examples are used to assess the 
accuracy, stability, and robustness of the proposed schemes in higher-dimensional 
and coupled settings.

\subsubsection{High-Dimensional Allen--Cahn Equation}

We consider two complementary twenty-dimensional Allen--Cahn problems.
The first problem admits a manufactured exact solution and is designed to
assess both approximation accuracy and robustness with respect to temporal
stiffness. The second problem is an unforced periodic Allen--Cahn equation
for which no closed-form transient solution is available; in this case, we
examine the evolution of the free energy to assess the dissipative behavior
of the numerical methods.

For the first test, we consider the forced Allen--Cahn equation
\begin{equation}
\label{eq:highdim_ac_forced}
    \partial_t u
    =
    \nu \Delta u
    -
    \frac{u^3-u}{\epsilon^2}
    +
    f_{\mathrm{AC}}(x,t),
    \qquad
    x\in\Omega=[-\pi,\pi]^d,
    \quad t>0,
\end{equation}
subject to periodic boundary conditions in every spatial direction. 
To construct a genuinely high-dimensional periodic solution, we define
\begin{equation*}
    \phi_d(x)
    =
    \frac{1}{\sqrt{2d}}
    \sum_{i=1}^{d}
    \left[
        \cos\bigl(q(x_i+x_{i+1})\bigr)
        +
        \cos\bigl(q(x_i-x_{i+1})\bigr)
    \right],
    \qquad
    x_{d+1}=x_1.
\end{equation*}
The manufactured exact solution of \eqref{eq:highdim_ac_forced} 
is chosen as
\begin{equation}
\label{eq:highdim_ac_exact}
    u_{\mathrm{ref}}(x,t)
    =
    1+w(x,t),
    \qquad
    w(x,t)
    =
    \delta e^{-\lambda t}\phi_d(x),
\end{equation}
where $\lambda=2\nu q^2+\frac{2}{\epsilon^2}$. 
Substitution of \eqref{eq:highdim_ac_exact} into
\eqref{eq:highdim_ac_forced} gives the forcing term
\begin{equation*}
    f_{\mathrm{AC}}(x,t)
    =
    \frac{3w(x,t)^2+w(x,t)^3}{\epsilon^2}.
\end{equation*}
Thus, the forcing compensates only for the nonlinear remainder, while the
stiff linear relaxation with characteristic rate
$2/\epsilon^2$ is retained in the dynamics.

We take $d=20, \, \nu=1, \, \epsilon=0.02, \, \delta=0.05, \, q=1$.
The time step and terminal time are $\Delta t=10^{-3}, \, T=1$. 
The periodic neural networks, Monte Carlo integration,
independent test points, and optimization budgets are specified
in \ref{app:highdim_settings}. The same numerical
configuration is used by all four methods. 
Table~\ref{tab:highdim_ac_accuracy} reports the
final-time relative $L^2$ errors evaluated against
the exact solution \eqref{eq:highdim_ac_exact}. As discussed in 
the preceding stability analysis, explicit time-stepping methods are 
subject to severe restrictions in stiff regimes. For $\epsilon=0.02$, 
the characteristic decay rate $\lambda=5002$ gives $\lambda\Delta t=5.002$, 
exceeding the classical forward-Euler stability limit and indicating 
the severe temporal stability challenge posed by the chosen parameters 
Consistent with this analysis, both EDNN and E-ONTS break down with 
NaN/Inf values, whereas I-ONTS and PC-ONTS remain 
numerically stable, achieving relative errors of $4.78\times10^{-2}$ 
and $5.11\times10^{-2}$, respectively. These results further demonstrate 
the improved robustness of the implicit Onsager formulations 
for high-dimensional stiff problems.

\begin{table}[H]
\centering
\caption{
Final-time relative $L^2$ errors for the twenty-dimensional
stiff forced Allen--Cahn equation with $\epsilon=0.02$
and $\Delta t=10^{-3}$. ``NaN/Inf'' indicates numerical breakdown during time evolution.
}
\label{tab:highdim_ac_accuracy}
\begin{tabular}{ccccc}
\toprule
Method & EDNN & E-ONTS & I-ONTS & PC-ONTS \\
\midrule
Relative $L^2$ error
& NaN/Inf
& NaN/Inf
& $\mathbf{4.7834\times10^{-2}}$
& $5.1146\times10^{-2}$ \\
\bottomrule
\end{tabular}
\end{table}

We next consider an unforced periodic Allen--Cahn equation,
\begin{equation}
\label{eq:highdim_ac_periodic}
    \partial_tu
    =
    \nu\Delta u
    -
    \frac{u^3-u}{\epsilon^2},
    \qquad
    x\in[-\pi,\pi]^d,
\end{equation}
subject to periodic boundary conditions. In this experiment, $d=20, \, \nu=1, \, \epsilon=0.2$ 
and the initial condition of \eqref{eq:highdim_ac_periodic} is chosen as
\begin{equation*}
    u_0(x)
    =
    \frac{A}{\sqrt d}
    \sum_{i=1}^{d}\cos(x_i),
    \qquad
    A=0.2.
\end{equation*}
All other experimental settings are the same 
as those in the preceding Allen--Cahn example 
with an exact solution.

The corresponding Allen--Cahn free energy is
\begin{equation}
\label{eq:highdim_ac_energy}
    \mathcal E(u)
    =
    \int_{\Omega}
    \left[
        \frac{\nu}{2}|\nabla u|^2
        +
        \frac{1}{4\epsilon^2}(u^2-1)^2
    \right]\,dx.
\end{equation}
Because the volume
$|\Omega|=(2\pi)^d$ becomes extremely large for $d=20$, we report the
mean free-energy density
\begin{equation}
\label{eq:highdim_ac_mean_energy}
    \overline{\mathcal E}(u)
    =
    \frac{\mathcal E(u)}{|\Omega|}.
\end{equation}
This normalization does not affect the energy-dissipation property.
For the Allen--Cahn gradient flow,
\begin{equation*}
    \frac{d}{dt}\mathcal E(u(t))
    =
    -
    \int_{\Omega}
    |\partial_tu|^2\,dx
    \leq0,
\end{equation*}
and hence $\overline{\mathcal E}(u)$ should also be nonincreasing.

\begin{figure}[H]
    \centering
    \includegraphics[width=0.65\textwidth]
    {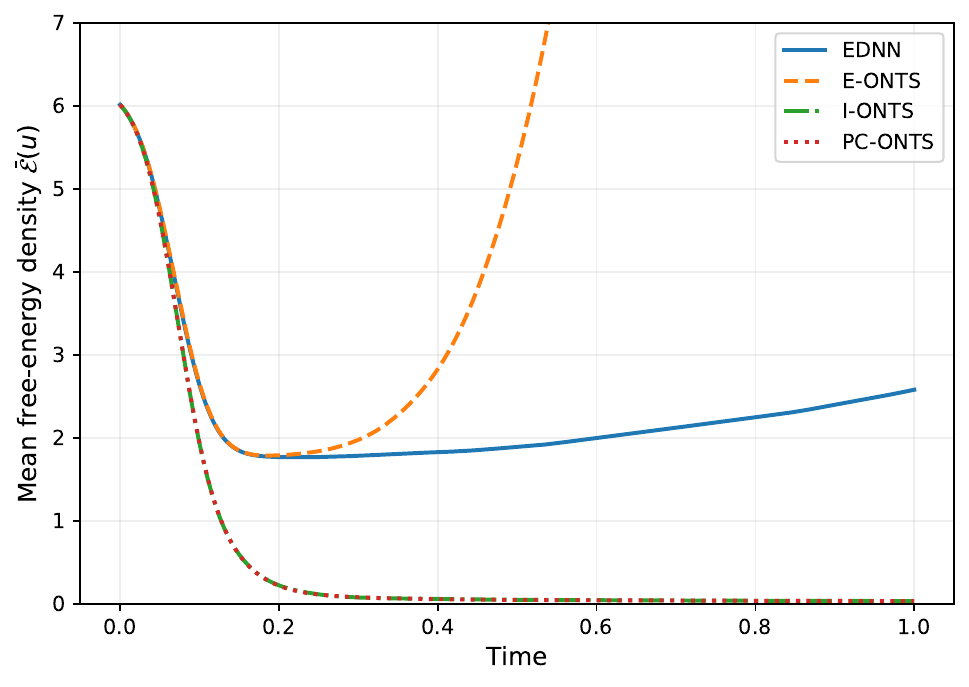}
    \caption{
    Evolution of the mean free-energy density for the
    twenty-dimensional periodic Allen--Cahn equation with
    $\Delta t=10^{-3}$.
    The vertical axis is restricted to $[0,7]$ for readability.
    }
    \label{fig:highdim_ac_energy}
\end{figure}

Since no closed-form transient solution is available, we focus 
on the free-energy evolution to assess whether the numerical methods 
preserve the expected dissipative behavior. Figure~\ref{fig:highdim_ac_energy} 
reveals a clear difference between the explicit and implicit temporal treatments. 
Although all four methods initially capture the rapid decrease in free energy, 
the energy curves of EDNN and E-ONTS subsequently turn upward, exhibiting 
significant nonphysical energy growth. In contrast, I-ONTS and PC-ONTS 
maintain monotonically decreasing mean free-energy densities throughout 
the simulation. This behavior is consistent with the discrete energy-dissipation 
property established in Proposition~\ref{prop:implicit-nonconserved-energy}, 
which is maintained by PC-ONTS through its safeguarded implicit correction.
These results demonstrate the robustness of 
the implicit Onsager formulations in preserving energy dissipation even in 
twenty-dimensional problems.

Together, the two high-dimensional Allen--Cahn experiments demonstrate 
the effectiveness of the implicit Onsager formulations in handling stiff 
dynamics and preserving energy dissipation in high-dimensional nonlinear gradient flows.
The manufactured-solution test confirms their numerical stability and accuracy 
in twenty dimensions, while the periodic experiment demonstrates their ability 
to preserve the expected energy-dissipation behavior without an exact transient 
solution. These findings further support the applicability of the proposed implicit 
schemes to high-dimensional problems.

\subsubsection{High-Dimensional Fokker--Planck Equation}

We next consider a $d$-dimensional periodic
Fokker--Planck equation with anisotropic time scales,
\begin{equation}
\label{eq:highdim_fp}
    \partial_t \rho
    =
    \sum_{i=1}^{d}
    \lambda_i
    \partial_{x_i}
    \left(
        \partial_{x_i}\rho
        +
        V'(x_i)\rho
    \right),
    \qquad
    \bm{x}\in\mathbb T^d=[0,1)^d,
\end{equation}
where
\begin{equation*}
    V(x)=\alpha\left(1-\cos(2\pi x)\right),
    \qquad
    \lambda_i=K^{\frac{i-1}{d-1}}.
\end{equation*}
The coefficients $\lambda_i$ range from $1$ to $K$,
introducing characteristic relative time scales
ranging from $1$ to $K^{-1}$. 
This choice introduces substantially 
different relaxation rates across spatial directions, 
providing a test of the methods under anisotropic dynamics.
The corresponding equilibrium density is
\begin{equation*}
    \rho_\infty(\bm{x})
    =
    \prod_{i=1}^{d}
    \frac{\exp\left(\alpha\cos(2\pi x_i)\right)}
         {I_0(\alpha)},
\end{equation*}
where $I_0$ denotes the modified Bessel function
of the first kind of order zero.
Introducing the density ratio
$u=\rho/\rho_\infty$ to \eqref{eq:highdim_fp} and using
$\partial_{x_i}\rho_\infty=-V'(x_i)\rho_\infty$, 
we obtain
\begin{equation}
\label{eq:highdim_fp_ratio}
    \partial_t u
    =
    \sum_{i=1}^{d}
    \lambda_i
    \left(
        \partial_{x_i x_i}u
        -
        V'(x_i)\partial_{x_i}u
    \right).
\end{equation}

For the Onsager formulations, we introduce the rescaled
flux variable
\begin{equation}
\label{eq:highdim_fp_flux}
    v_i=-\partial_{x_i}u,
    \qquad i=1,\ldots,d,
\end{equation}
which is related to the physical flux through
$j_i=\lambda_i\rho_\infty v_i$. 
The corresponding evolution equation can be written as
\begin{equation*}
    \partial_t u
    +
    \sum_{i=1}^{d}
    \lambda_i
    \left(
        \partial_{x_i}v_i
        -
        V'(x_i)v_i
    \right)
    =0.
\end{equation*}
This rescaled formulation is equivalent to the conservative 
flux representation of the physical density $\rho=u\rho_\infty$ 
and preserves the weighted mass
$\int_{\mathbb T^d}u\rho_\infty\,\mathrm{d}\bm{x}$.
We monitor the weighted quadratic energy
\begin{equation*}
    \mathcal E(u)
    =
    \frac{1}{2}
    \int_{\mathbb T^d}
    (u-1)^2\rho_\infty\,\mathrm{d}\bm{x},
\end{equation*}
which satisfies the continuous dissipation law
\begin{equation*}
    \frac{d}{dt}\mathcal E(u)
    =
    -\sum_{i=1}^{d}
    \lambda_i
    \int_{\mathbb T^d}
    \rho_\infty
    |\partial_{x_i}u|^2\,\mathrm{d}\bm{x}
    \leq 0.
\end{equation*}
The initial density ratio is prescribed as
\begin{equation*}
    u_0(\bm{x})
    =
    1
    +
    \frac{\eta}{d}
    \sum_{i=1}^{d}
    \left[
        \cos(2\pi x_i)-c_\alpha
    \right],
    \qquad
    c_\alpha=\frac{I_1(\alpha)}{I_0(\alpha)},
\end{equation*}
where $I_1$ denotes the modified Bessel function
of the first kind of order one.
The constant $c_\alpha$ ensures that the perturbation 
has zero mean with respect to the equilibrium measure, 
so that the initial weighted mass is preserved.
We set $d=20$, $\alpha=1$, $\eta=0.5$,
$K=10$, $\Delta t=10^{-3}$, and $T=1$. 
Using the separable structure, a high-accuracy reference solution
is assembled from a one-dimensional auxiliary problem:
\begin{equation}
\label{eq:highdim_fp_reference}
    u_{\mathrm{ref}}(\bm{x},t)
    =1+\frac{\eta}{d}\sum_{i=1}^{d}h(x_i,\lambda_i t),
\end{equation}
where $h$ solves the problem given in
\ref{app:highdim_settings}. 
Spatial integrals are evaluated by Monte Carlo sampling from
$\rho_\infty$, with an independent sample set for the reported
space-time relative $L^2(\rho_\infty)$ errors. Full details of
sampling, network architectures, and optimization are provided in
\ref{app:highdim_settings}.

Table~\ref{tab:highdim_fp_accuracy} shows a substantial
accuracy advantage of the implicit Onsager formulations
for the high-dimensional anisotropic Fokker--Planck equation.
The space-time relative $L^2(\rho_\infty)$ errors of EDNN
and E-ONTS are $9.48\times10^{-2}$ and
$1.79\times10^{-1}$, respectively, whereas both I-ONTS
and PC-ONTS achieve an error of approximately
$5.19\times10^{-3}$. 
In this experiment, the explicit predictor in PC-ONTS
is rejected by the Rayleighian-based safeguard
\eqref{eq:Rayleighian-based_safeguard} at every time step.
Consequently, the previous parameter vector is used
to initialize the implicit correction, as in I-ONTS.
This explains the nearly identical results obtained
by the two methods and indicates that the explicit
predictor provides little benefit under the present
anisotropic time scales.

\begin{table}[H]
\centering
\caption{
Space-time relative $L^2(\rho_\infty)$ errors for the
twenty-dimensional periodic Fokker--Planck equation.
}
\label{tab:highdim_fp_accuracy}
\begin{tabular}{ccccc}
\toprule
Method & EDNN & E-ONTS & I-ONTS & PC-ONTS \\
\midrule
Relative $L^2(\rho_\infty)$ error
& $9.4821\times10^{-2}$
& $1.7922\times10^{-1}$
& $\mathbf{5.1856\times10^{-3}}$
& $\mathbf{5.1856\times10^{-3}}$ \\
\bottomrule
\end{tabular}
\end{table}

\begin{figure}[H]
    \centering
    \begin{subfigure}{0.48\textwidth}
        \centering
        \includegraphics[width=\textwidth]
        {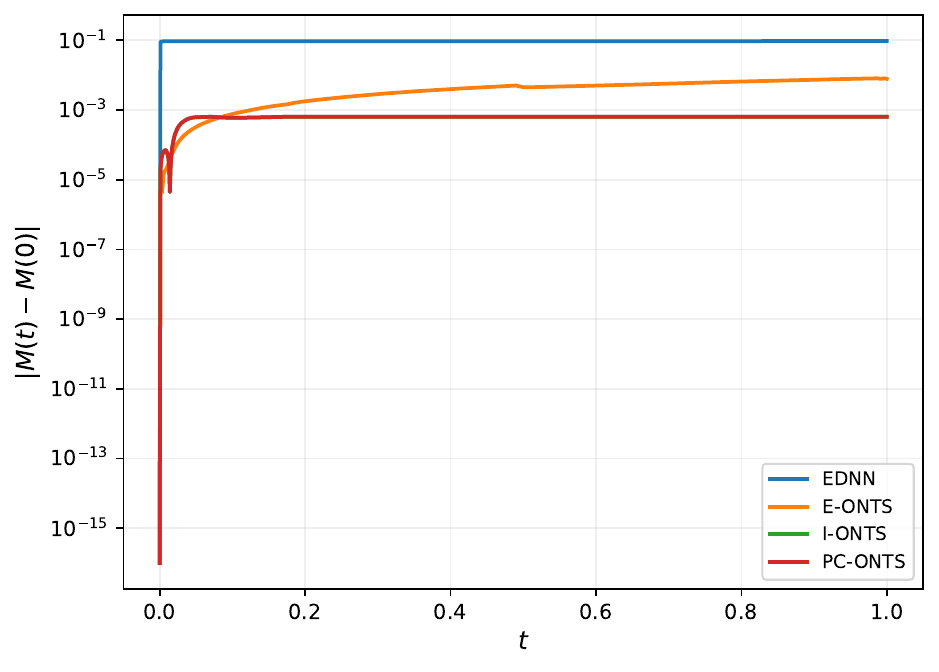}
        \caption{Absolute mass error.}
        \label{fig:highdim_fp_mass_error}
    \end{subfigure}
    \hfill
    \begin{subfigure}{0.48\textwidth}
        \centering
        \includegraphics[width=\textwidth]
        {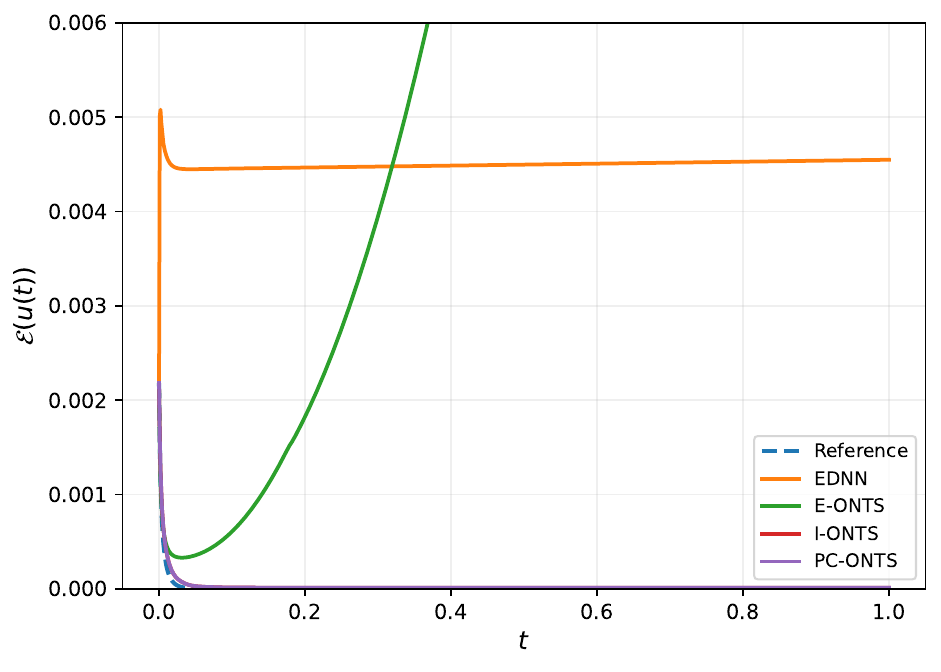}
        \caption{Weighted free-energy evolution.}
        \label{fig:highdim_fp_energy}
    \end{subfigure}

    \caption{
    Mass conservation and weighted free-energy evolution for the
    twenty-dimensional periodic Fokker--Planck equation with
    $K=10$ and $\Delta t=10^{-3}$.
    }
    \label{fig:highdim_fp_mass_energy}
\end{figure}

Figure~\ref{fig:highdim_fp_mass_energy} further compares
the mass conservation and energy dissipation properties.
The physical mass is defined by
\begin{equation}
\label{eq:highdim_fp_mass}
    M(t)
    =
    \int_{\mathbb T^d}\rho(\bm{x},t)\,\mathrm{d}\bm{x}
    =
    \int_{\mathbb T^d}
    u(\bm{x},t)\rho_\infty(\bm{x})\,\mathrm{d}\bm{x},
\end{equation}
and the mass error is measured by
\begin{equation}
\label{eq:highdim_fp_mass_error}
    \mathrm{Err}_M(t)
    =
    |M(t)-M(0)|.
\end{equation}
As shown in Figure~\ref{fig:highdim_fp_mass_error},
I-ONTS and PC-ONTS exhibit substantially smaller
numerically estimated mass errors than EDNN and E-ONTS.
Although mass conservation is inherent in the
flux-based ONTS formulations, the computed mass
may exhibit deviations due to Monte Carlo
integration errors, particularly when the numerical
solution becomes inaccurate or unstable.
Figure~\ref{fig:highdim_fp_energy} shows the evolution
of the weighted free energy. I-ONTS and PC-ONTS
closely follow the reference energy decay throughout
the simulation, consistent with the discrete
energy-dissipation property of the implicit
Onsager formulation.
In contrast, EDNN deviates significantly from
the reference evolution, while E-ONTS exhibits
rapid energy growth after an initial decrease
and eventually exceeds the plotting range.
Together, these results demonstrate the improved
accuracy and structure-preserving performance
of the implicit Onsager methods for high-dimensional
anisotropic dynamics.

\subsubsection{Coupled Multi-Species Poisson--Nernst--Planck System}
We conclude the numerical study with a coupled three-species
Poisson--Nernst--Planck (PNP) system. Unlike the preceding single-field
examples, this problem involves several conserved concentrations coupled
through a common electrostatic potential. It therefore provides a test of
the proposed schemes for nonlinear, multi-component gradient flows.
On the periodic domain $\Omega=\mathbb T^2=[0,2\pi)^2$, we consider

\begin{equation}
\label{eq:pnp3_system}
\begin{cases}
    \partial_t c_i+\nabla\cdot \bm{j}_i =0,&\\
    \bm{j}_i=-D_i\left(\nabla c_i+z_i c_i\nabla\phi\right),
    &\qquad i=1,2,3,\\
    -\varepsilon_p^2\Delta\phi =\sum_{i=1}^3 z_i c_i,&
\end{cases}
\end{equation}
where $c_i$ and $\bm{j}_i$ denote the concentration and flux of the $i$-th
species, respectively, $D_i>0$ is its diffusivity, $z_i$ is its valence,
and $\phi$ is the electrostatic potential. Periodic boundary conditions
are imposed on all fields, and the potential is fixed by the zero-mean
gauge $\int_\Omega\phi\,\mathrm{d}\bm x=0$. We take
\begin{equation*}
    (z_1,z_2,z_3)=(1,-1,2),\qquad
    (D_1,D_2,D_3)=(1,1,1),\qquad
    \varepsilon_p=0.05.
\end{equation*}
The small Debye-length parameter introduces strong
electrostatic coupling. Linearization around the
electroneutral equilibrium gives a characteristic
charge-relaxation rate of order
$\varepsilon_p^{-2}\sum_{i=1}^3 D_i z_i^2c_{i,\infty}=3200$, 
leading to a restrictive time-step condition
for explicit discretizations.

The homogeneous equilibrium of \eqref{eq:pnp3_system} is given by
\begin{equation*}
    \bm c_\infty=(1,3,1)^T,\qquad \phi_\infty=0,
\end{equation*}
which satisfies the electroneutrality condition
$\sum_{i=1}^3z_i c_{i,\infty}=0$.
We prescribe the initial concentrations as
\begin{equation}
\label{eq:pnp3_initial}
\begin{aligned}
    c_1(x,y,0)&=1+a\cos x\cos y,\\
    c_2(x,y,0)&=3+a\cos(x+y),\\
    c_3(x,y,0)&=1+a\sin x\sin y,
\end{aligned}
\qquad a=0.2.
\end{equation}
Each perturbation has zero spatial mean, so the individual species
masses agree with their equilibrium values and the total charge is zero.
Nevertheless, the initial charge density is spatially nonuniform,
producing a nontrivial electrostatic field. In fact, the initial
potential is explicitly given by
\begin{equation*}
    \phi(x,y,0)
    =\frac{3a}{4\varepsilon_p^2}
    \bigl[\cos(x-y)-\cos(x+y)\bigr].
\end{equation*}
For positive concentrations, the PNP system 
admits the following relative free-energy functional
\begin{equation}
\label{eq:pnp3_energy}
    \mathcal E(\bm c)
    =\int_\Omega
    \left[
        \sum_{i=1}^3
        \left(
            c_i\log\frac{c_i}{c_{i,\infty}}
            -c_i+c_{i,\infty}
        \right)
        +\frac{\varepsilon_p^2}{2}|\nabla\phi|^2
    \right]\mathrm{d}\bm x.
\end{equation}
The corresponding electrochemical potentials are
$\mu_i=\log(c_i/c_{i,\infty})+z_i\phi$, and hence
$\bm{j}_i=-D_i c_i\nabla\mu_i$.
The dissipation potential takes the form
\begin{equation*}
    \Phi(\bm c;\bm{j})
    =\frac12\sum_{i=1}^3
    \int_\Omega\frac{|\bm {j}_i|^2}{D_i c_i}\,\mathrm{d}\bm x.
\end{equation*}
Consequently, smooth positive solutions satisfy
\begin{equation*}
    \frac{d}{dt}\mathcal E(\bm c(t))
    =-\sum_{i=1}^3\int_\Omega
        D_i c_i|\nabla\mu_i|^2\,\mathrm{d}\bm x
    =-2\Phi(\bm c;\bm{j})\leq0.
\end{equation*}
The equilibrium energy is $\mathcal E(\bm c_\infty)=0$.
In the ONTS schemes, a separate neural flux is introduced
for each species, and the concentrations are advanced through
\begin{equation}
\label{eq:pnp3_conservative_update}
    c_i^{k+1}
    =
    c_i^k-\Delta t\,\nabla\cdot
    \bm{j}_{i,\bm\theta^\star},
    \qquad i=1,2,3,
\end{equation}
where $\bm\theta^\star=\bm\theta^k$ for E-ONTS and
$\bm\theta^\star=\bm\theta^{k+1}$ for I-ONTS and PC-ONTS.
For I-ONTS and PC-ONTS, the implicit correction approximately
minimizes the discrete Onsager objective
\begin{equation*}
    \mathcal J_k(\bm\theta)
    =\mathcal E\!\left(
        \bm c^k-\Delta t\,\nabla\cdot\bm{j}_{\bm\theta}
    \right)
    +\frac{\Delta t}{2}\sum_{i=1}^3
    \int_\Omega
        \frac{|\bm{j}_{i,\bm\theta}|^2}{D_i c_i^k}\,\mathrm{d}\bm x,
\end{equation*}
with the mobility $D_i c_i^k$ evaluated at the previous time level.
In PC-ONTS, the explicit predictor is screened by the Rayleighian-based
safeguard \eqref{eq:Rayleighian-based_safeguard} before the implicit
correction. The zero-flux candidate provides the energy-descent safeguard
described in \eqref{eq:implicit-safeguarded-selection}.

We set $\Delta t=10^{-3}$ and $T=1$. All neural solvers use a
$64\times64$ periodic collocation grid. The reference solution is
obtained using a Fourier pseudospectral method on a separate
$96\times96$ grid. Network architectures, Poisson zero-mode
handling, positivity treatment, and the space-time error measure
\eqref{eq:pnp3_spacetime_error} are given in
\ref{app:pnp_settings}. 
Table~\ref{tab:pnp3_accuracy} shows that the implicit Onsager
formulations substantially improve the solution accuracy for this
strongly coupled system. EDNN and E-ONTS yield space-time relative
$L^2$ errors of $1.0046\times10^{-1}$ and $2.0035\times10^{-1}$,
respectively, whereas I-ONTS reduces the error to
$3.8447\times10^{-2}$. PC-ONTS achieves the smallest error,
$2.8368\times10^{-2}$, representing a further reduction of
approximately $26\%$ compared with I-ONTS.

\begin{table}[H]
\centering
\caption{Space-time relative $L^2$ errors for the coupled three-species
Poisson--Nernst--Planck system.}
\label{tab:pnp3_accuracy}
\small
\setlength{\tabcolsep}{4pt}
\begin{tabular}{@{}lcccc@{}}
\toprule
Method & EDNN & E-ONTS & I-ONTS & PC-ONTS \\
\midrule
Relative $L^2$ error
& $1.0046\times10^{-1}$
& $2.0035\times10^{-1}$
& $3.8447\times10^{-2}$
& $\mathbf{2.8368\times10^{-2}}$ \\
\bottomrule
\end{tabular}
\end{table}

\begin{figure}[H]
    \centering
    \begin{subfigure}[t]{0.48\linewidth}
        \centering
        \includegraphics[width=\linewidth]{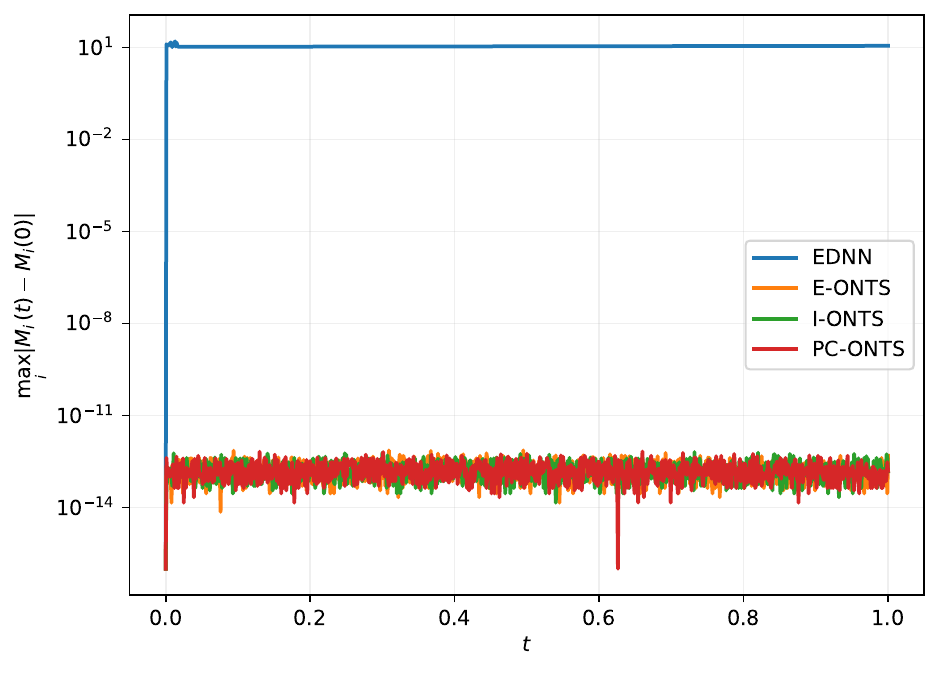}
        \caption{Maximum species-wise mass error.}
        \label{fig:pnp3_mass_error}
    \end{subfigure}
    \hfill
    \begin{subfigure}[t]{0.48\linewidth}
        \centering
        \includegraphics[width=\linewidth]{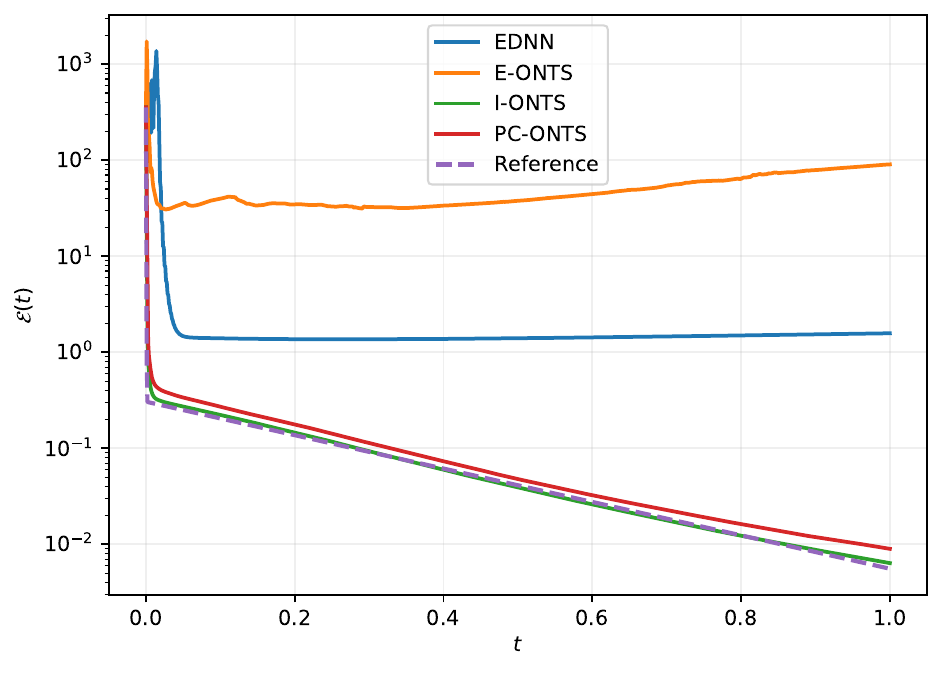}
        \caption{Relative free-energy evolution.}
        \label{fig:pnp3_energy}
    \end{subfigure}
    \caption{Comparison of mass conservation and free-energy
    dissipation for the coupled three-species
    Poisson--Nernst--Planck system:
    (a) maximum species-wise mass error;
    (b) evolution of the relative free energy compared
    with the Fourier pseudospectral reference solution.}
    \label{fig:pnp3_mass_energy}
\end{figure}

To assess the conservation properties, we monitor the mass of each
species and its deviation from the initial value,
\begin{equation*}
    M_i(t)=\int_\Omega c_i(\bm x,t)\,\mathrm{d}\bm x,
    \qquad
    \mathrm{Err}_{M,i}(t)=|M_i(t)-M_i(0)|,
    \quad i=1,2,3.
\end{equation*}
Figure~\ref{fig:pnp3_mass_energy}(a) displays the maximum
species-wise mass error, $\max_i\mathrm{Err}_{M,i}(t)$.
All three ONTS methods, including the explicit E-ONTS scheme,
maintain mass errors at approximately the $10^{-13}$ level,
whereas EDNN exhibits a much larger mass drift, reaching
the order of $10^1$.
This contrast reflects the flux-based update
\eqref{eq:pnp3_conservative_update}, which conserves the mass
of each species under periodic boundary conditions,
independently of the accuracy of the learned fluxes.
The small nonzero errors of the ONTS methods are attributable
to spatial quadrature and floating-point effects. 
Figure~\ref{fig:pnp3_mass_energy}(b) compares the evolution
of the relative free energy \eqref{eq:pnp3_energy} with 
that of the reference solution.
I-ONTS and PC-ONTS exhibit monotonically decreasing free
energies that closely follow the reference dissipation trend.
This behavior is consistent with the discrete energy-dissipation
property established in
Proposition~\ref{prop:implicit-conserved-structure},
under the corresponding descent acceptance criterion.
In contrast, EDNN captures the initial rapid energy decrease
but subsequently stagnates at an energy level well above
the reference.
E-ONTS also exhibits an initial decrease, followed by
pronounced energy growth as the simulation proceeds.
Thus, although E-ONTS retains the species-wise mass
conservation property, its explicit temporal treatment
does not maintain the expected dissipative behavior
in this stiff electrostatically coupled problem.

This example extends the numerical study beyond scalar gradient
flows to a coupled system with three individually conserved
concentrations, state-dependent mobilities, and nonlinear
electrostatic interactions.
Together, the accuracy, mass-conservation, and
energy-dissipation results demonstrate the effectiveness
of the implicit Onsager schemes for multi-component
gradient-flow problems, while highlighting that
conservation by itself does not guarantee stable
energy dissipation.

\section{Conclusions}

In this work, we developed a unified Onsager-based neural
time-stepping (ONTS) framework for conserved and nonconserved
gradient-flow PDEs. Rather than directly parameterizing the
evolving state, the proposed framework employs spatial neural
networks to approximate the Onsager process variables, namely, 
the rates in nonconserved systems and the fluxes in conserved systems. 
Within this framework, we introduced three complementary
schemes: an explicit least-squares scheme (E-ONTS), an implicit
variational scheme (I-ONTS), and a predictor--corrector scheme
(PC-ONTS) that combines explicit prediction with safeguarded
implicit correction. The proposed framework incorporates important structural
properties of the underlying gradient flows. For conserved
systems, the flux-based update preserves mass under periodic
or no-flux boundary conditions, independently of neural-network
optimization accuracy. For I-ONTS and PC-ONTS, we rigorously 
established unconditional discrete energy stability under 
a suitable descent criterion, without requiring global 
minimization of the neural variational objective. 
These properties provide a theoretical foundation for
structure-preserving neural time integration.

Numerical experiments on heat, Allen--Cahn, Cahn--Hilliard,
Fokker--Planck, and coupled Poisson--Nernst--Planck equations
demonstrated the accuracy, robustness and structure-preserving properties
of the proposed methods. While E-ONTS
remains competitive for relatively smooth problems, the
implicit schemes exhibit improved robustness in stiff,
high-dimensional, and coupled nonlinear systems.
The predictor--corrector strategy can further improve
accuracy or reduce the optimization effort in suitable
settings, although its effectiveness depends on the
reliability of the explicit predictor.
Future work will explore higher-order structure-preserving
time discretizations, positivity-preserving extensions
for systems involving density or concentration variables,
and more efficient neural optimization strategies for
large-scale problems. Another promising direction is
to extend the proposed framework to structure-preserving
neural operator learning, with the aim of retaining
mass conservation and discrete energy stability across
different initial conditions and model parameters.

\section*{Acknowledgments}

This work was supported in part by the National Natural
Science Foundation of China
(Grant Nos.~12288201, 12371415, 12461160275, and 12625117),
the Science Challenge Project (Grant No.~TZ2025006),
the Strategic Priority Research Program of the Chinese
Academy of Sciences (Grant No.~XDA0480504),
the NSFC/RGC Joint Research Scheme
(Grant No.~N\_PolyU5145/24),
the Research Institute for Smart Energy of
The Hong Kong Polytechnic University
(Project No.~4-CDDQ),
and the Beijing Natural Science Foundation
(Grant No.~Z240001).

\bibliographystyle{elsarticle-num}
\bibliography{ref}

@article{raissi2019physics,
  title={Physics-informed neural networks: A deep learning framework for solving forward and inverse problems involving nonlinear partial differential equations},
  author={Raissi, Maziar and Perdikaris, Paris and Karniadakis, George E},
  journal={Journal of Computational Physics},
  volume={378},
  pages={686--707},
  year={2019},
  publisher={Elsevier}
}

@article{pang2019fpinns,
  title={{FPINNs}: Fractional physics-informed neural networks},
  author={Pang, Guofei and Lu, Lu and Karniadakis, George Em},
  journal={SIAM Journal on Scientific Computing},
  volume={41},
  number={4},
  pages={A2603--A2626},
  year={2019},
  publisher={SIAM}
}

@article{jagtap2020extended,
  title={Extended physics-informed neural networks {(XPINNs)}: A generalized space-time domain decomposition based deep learning framework for nonlinear partial differential equations},
  author={Jagtap, Ameya D and Karniadakis, George Em},
  journal={Communications in Computational Physics},
  volume={28},
  number={5},
  year={2020},
  publisher={Brown Univ., Providence, RI (United States)}
}

@article{yu2018deep,
  title={The deep {Ritz} method: A deep learning-based numerical algorithm for solving variational problems},
  author={Yu, Bing and others},
  journal={Communications in Mathematics and Statistics},
  volume={6},
  number={1},
  pages={1--12},
  year={2018},
  publisher={Springer}
}

@article{liao2021,
  author = {Liao, Yulei and Ming, Pingbing},
  title = {Deep {Nitsche} method: Deep {Ritz} method with essential boundary conditions},
  journal = {Communications in Computational Physics},
  year = {2021},
  volume = {29},
  number = {5},
  pages = {1365-1384},
}

@article{sirignano2018dgm,
  title={{DGM}: A deep learning algorithm for solving partial differential equations},
  author={Sirignano, Justin and Spiliopoulos, Konstantinos},
  journal={Journal of Computational Physics},
  volume={375},
  pages={1339--1364},
  year={2018},
  publisher={Elsevier}
}

@article{zang2020weak,
  title={Weak adversarial networks for high-dimensional partial differential equations},
  author={Zang, Yaohua and Bao, Gang and Ye, Xiaojing and Zhou, Haomin},
  journal={Journal of Computational Physics},
  volume={411},
  pages={109409},
  year={2020},
  publisher={Elsevier}
}

@article{kharazmi2021hp,
  title={{hp-VPINNs}: Variational physics-informed neural networks with domain decomposition},
  author={Kharazmi, Ehsan and Zhang, Zhongqiang and Karniadakis, George Em},
  journal={Computer Methods in Applied Mechanics and Engineering},
  volume={374},
  pages={113547},
  year={2021},
  publisher={Elsevier}
}

@article{han2017deep,
  title={Deep learning-based numerical methods for high-dimensional parabolic partial differential equations and backward stochastic differential equations},
  author={Han, Jiequn and Jentzen, Arnulf and others},
  journal={Communications in Mathematics and Statistics},
  volume={5},
  number={4},
  pages={349--380},
  year={2017},
  publisher={Springer}
}

@article{chen2022bridging,
  title={Bridging traditional and machine learning-based algorithms for solving {PDEs}: The random feature method},
  author={Chen, Jingrun and Chi, Xurong and Yang, Zhouwang and others},
  journal={Journal of Machine Learning},
  volume={1},
  number={3},
  pages={268--298},
  year={2022}
}

@article{lu2021learning,
  title={Learning nonlinear operators via {DeepONet} based on the universal approximation theorem of operators},
  author={Lu, Lu and Jin, Pengzhan and Pang, Guofei and Zhang, Zhongqiang and Karniadakis, George Em},
  journal={Nature Machine Intelligence},
  volume={3},
  number={3},
  pages={218--229},
  year={2021},
  publisher={Nature Publishing Group UK, London}
}

@article{wang2021learning,
  title={Learning the solution operator of parametric partial differential equations with physics-informed {DeepONets}},
  author={Wang, Sifan and Wang, Hanwen and Perdikaris, Paris},
  journal={Science Advances},
  volume={7},
  number={40},
  pages={eabi8605},
  year={2021},
  publisher={American Association for the Advancement of Science}
}

@article{wang2022improved,
  title={Improved architectures and training algorithms for deep operator networks},
  author={Wang, Sifan and Wang, Hanwen and Perdikaris, Paris},
  journal={Journal of Scientific Computing},
  volume={92},
  number={2},
  pages={35},
  year={2022},
  publisher={Springer}
}

@inproceedings{li2021fourier,
  title={Fourier Neural Operator for Parametric Partial Differential Equations},
  author={Li, Zongyi and Kovachki, Nikola and
          Azizzadenesheli, Kamyar and Liu, Burigede and
          Bhattacharya, Kaushik and Stuart, Andrew and
          Anandkumar, Anima},
  booktitle={International Conference on Learning Representations},
  year={2021},
}

@article{kovachki2023neural,
  title={Neural operator: Learning maps between function spaces with applications to {PDEs}},
  author={Kovachki, Nikola and Li, Zongyi and Liu, Burigede and Azizzadenesheli, Kamyar and Bhattacharya, Kaushik and Stuart, Andrew and Anandkumar, Anima},
  journal={Journal of Machine Learning Research},
  volume={24},
  number={89},
  pages={1--97},
  year={2023}
}

@article{gupta2021multiwavelet,
  title={Multiwavelet-based operator learning for differential equations},
  author={Gupta, Gaurav and Xiao, Xiongye and Bogdan, Paul},
  journal={Advances in Neural Information Processing Systems},
  volume={34},
  pages={24048--24062},
  year={2021}
}

@article{li2024physics,
  title={Physics-informed neural operator for learning partial differential equations},
  author={Li, Zongyi and Zheng, Hongkai and Kovachki, Nikola and Jin, David and Chen, Haoxuan and Liu, Burigede and Azizzadenesheli, Kamyar and Anandkumar, Anima},
  journal={ACM/IMS Journal of Data Science},
  volume={1},
  number={3},
  pages={1--27},
  year={2024},
  publisher={ACM, New York, NY}
}

@article{lan2023dosnet,
  title={{DOSnet} as a non-black-box {PDE} solver: When deep learning meets operator splitting},
  author={Lan, Yuan and Li, Zhen and Sun, Jie and Xiang, Yang},
  journal={Journal of Computational Physics},
  volume={491},
  pages={112343},
  year={2023},
  publisher={Elsevier}
}

@article{huang2027pidosnet,
  title   = {{PI-DOSnet}: A physics-informed deep operator-splitting network for evolution partial differential equations},
  author  = {Huang, Jizu and Qian, Yue and Zhou, Tao},
  journal = {Journal of Computational Physics},
  volume  = {570},
  pages   = {115552},
  year    = {2027},
  publisher={Elsevier}
}

@article{jagtap2020conservative,
  title={Conservative physics-informed neural networks on discrete domains for conservation laws: Applications to forward and inverse problems},
  author={Jagtap, Ameya D and Kharazmi, Ehsan and Karniadakis, George Em},
  journal={Computer Methods in Applied Mechanics and Engineering},
  volume={365},
  pages={113028},
  year={2020},
  publisher={Elsevier}
}

@article{huang2026mass,
  title={Mass-preserving spatio-temporal adaptive {PINN} for {Cahn-Hilliard} equations with strong nonlinearities and singularities},
  author={Huang, Qiumei and Ma, Jiaxuan and Xu, Zhen},
  journal={Neural Networks},
  pages={109634},
  year={2026},
  publisher={Elsevier}
}

@article{badia2025compatible,
  title={Compatible finite element interpolated neural networks},
  author={Badia, Santiago and Li, Wei and Mart{\'\i}n, Alberto F},
  journal={Computer Methods in Applied Mechanics and Engineering},
  volume={439},
  pages={117889},
  year={2025},
  publisher={Elsevier}
}

@article{boon2025neural,
  title={Neural network solvers for parametrized elasticity problems that conserve linear and angular momentum},
  author={Boon, Wietse M and Franco, Nicola R and Fumagalli, Alessio},
  journal={Computer Methods in Applied Mechanics and Engineering},
  volume={437},
  pages={117759},
  year={2025},
  publisher={Elsevier}
}

@article{kutuk2025energy,
  title={Energy dissipation preserving physics informed neural network for {Allen--Cahn} equations},
  author={K{\"u}t{\"u}k, Mustafa and Y{\"u}cel, Hamdullah},
  journal={Journal of Computational Science},
  volume={87},
  pages={102577},
  year={2025},
  publisher={Elsevier}
}

@article{wu2023energy,
  title={Energy stable neural network for gradient flow equations},
  author={Wu, Yue and Jin, Tianyu and Chen, Chuqi and Fan, Ganghua and Lan, Yuan and Zhang, Luchan and Xiang, Yang},
  journal={arXiv preprint arXiv:2309.10002},
  year={2023}
}

@misc{huang2026edp,
  title={{EDP-PI-DOSnet}: An Energy Dissipation Preserving Neural Solver for Two-Phase Flow Problems},
  author={Huang, Jizu and Qian, Yue and Zhou, Tao},
  year={2026},
  eprint={2611.09430},
  archivePrefix={arXiv},
  primaryClass={math.NA},
}

@article{bamdad2026penco,
  title={Penco: A physics--energy--numerics--consistent operator for 3d phase field modeling},
  author={Bamdad, Mostafa and Eshaghi, Mohammad Sadegh and Anitescu, Cosmin and Valizadeh, Navid and Rabczuk, Timon},
  journal={Computer Methods in Applied Mechanics and Engineering},
  volume={454},
  pages={118862},
  year={2026},
  publisher={Elsevier}
}

@article{geng2024deep,
  title={A deep learning method for the dynamics of classic and conservative {Allen-Cahn} equations based on fully-discrete operators},
  author={Geng, Yuwei and Teng, Yuankai and Wang, Zhu and Ju, Lili},
  journal={Journal of Computational Physics},
  volume={496},
  pages={112589},
  year={2024},
  publisher={Elsevier}
}

@inproceedings{tanaka2025energy,
  title={Energy-consistent Neural Operators for {Hamiltonian} and Dissipative Partial Differential Equations.},
  author={Tanaka, Yusuke and Yaguchi, Takaharu and Iwata, Tomoharu and Ueda, Naonori},
  booktitle={AISTATS},
  pages={1882--1890},
  year={2025}
}

@article{chen2026structure,
  title={A Structure-Informed Framework for Solving Time-Dependent Partial Differential Equations With Neural Networks},
  author={Chen, Gaohang and Ju, Lili and Qiao, Zhonghua},
  journal={Numerical Methods for Partial Differential Equations},
  volume={42},
  number={4},
  pages={e70133},
  year={2026},
  publisher={Wiley Online Library}
}

@article{eyre1998unconditionally,
  title={Unconditionally gradient stable time marching the {Cahn-Hilliard} equation},
  author={Eyre, David J},
  journal={MRS Online Proceedings Library (OPL)},
  volume={529},
  pages={39},
  year={1998},
  publisher={Cambridge University Press}
}

@article{furihata2001stable,
  title={A stable and conservative finite difference scheme for the {Cahn-Hilliard} equation},
  author={Furihata, Daisuke},
  journal={Numerische Mathematik},
  volume={87},
  number={4},
  pages={675--699},
  year={2001},
  publisher={Springer}
}

@article{shen2018scalar,
  title={The scalar auxiliary variable {(SAV)} approach for gradient flows},
  author={Shen, Jie and Xu, Jie and Yang, Jiang},
  journal={Journal of Computational Physics},
  volume={353},
  pages={407--416},
  year={2018},
  publisher={Elsevier}
}

@article{hou2024energy,
  title={Energy-dissipative spectral renormalization exponential integrator method for gradient flow problems},
  author={Hou, Dianming and Ju, Lili and Qiao, Zhonghua},
  journal={SIAM Journal on Scientific Computing},
  volume={46},
  number={6},
  pages={A3477--A3502},
  year={2024},
  publisher={SIAM}
}

@article{chen2025onsager,
  title={The {Onsager} principle and structure preserving numerical schemes},
  author={Chen, Huangxin and Liu, Hailiang and Xu, Xianmin},
  journal={Journal of Computational Physics},
  volume={523},
  pages={113679},
  year={2025},
  publisher={Elsevier}
}

@article{li2026high,
  title={High-order unconditionally energy-stable decoupled method for phase field modeling of pitting corrosion with adaptive implementation},
  author={Li, Futuan and Li, Hongwei and Tang, Tao and Yang, Jiang},
  journal={Advances in Applied Mathematics and Mechanics},
  volume={18},
  number={5},
  pages={1355--1388},
  year={2026}
}

@article{zhang2026up,
  title={Up to Fourth-Order Structure-Preserving Diagonally Implicit {Runge-Kutta} Schemes for the {Allen-Cahn}-Type Equations},
  author={Zhang, Wei and Zhu, Xiaohan and Cai, Wenjun and Wang, Yushun},
  journal={Advances in Applied Mathematics and Mechanics},
  volume={18},
  number={6},
  pages={1914--1938},
  year={2026}
}

@article{huang2026efficient,
  title={Efficient bound-preserving numerical schemes for a phase-field model of tumor growth with extracellular matrix degradation},
  author={Huang, Qiumei and Qiao, Zhonghua and Wang, Cheng and Yang, Huiting},
  journal={Computer Methods in Applied Mechanics and Engineering},
  volume={450},
  pages={118621},
  year={2026},
  publisher={Elsevier}
}

@article{du2021evolutional,
  title={Evolutional deep neural network},
  author={Du, Yifan and Zaki, Tamer A},
  journal={Physical Review E},
  volume={104},
  number={4},
  pages={045303},
  year={2021},
  publisher={APS}
}

@article{kast2024positional,
  title={Positional embeddings for solving {PDEs} with evolutional deep neural networks},
  author={Kast, Mariella and Hesthaven, Jan S},
  journal={Journal of Computational Physics},
  volume={508},
  pages={112986},
  year={2024},
  publisher={Elsevier}
}

@article{bruna2024neural,
  title={Neural {Galerkin} schemes with active learning for high-dimensional evolution equations},
  author={Bruna, Joan and Peherstorfer, Benjamin and Vanden-Eijnden, Eric},
  journal={Journal of Computational Physics},
  volume={496},
  pages={112588},
  year={2024},
  publisher={Elsevier}
}

@article{georgoulis2023discrete,
  title={Discrete gradient flow approximations of high dimensional evolution partial differential equations via deep neural networks},
  author={Georgoulis, Emmanuil H and Loulakis, Michail and Tsiourvas, Asterios},
  journal={Communications in Nonlinear Science and Numerical Simulation},
  volume={117},
  pages={106893},
  year={2023},
  publisher={Elsevier}
}

@article{mattey2025gradient,
  title={Gradient flow based phase-field modeling using separable neural networks},
  author={Mattey, Revanth and Ghosh, Susanta},
  journal={Computer Methods in Applied Mechanics and Engineering},
  volume={439},
  pages={117897},
  year={2025},
  publisher={Elsevier}
}

@article{zhang2024energy,
  title={Energy-dissipative evolutionary deep operator neural networks},
  author={Zhang, Jiahao and Zhang, Shiheng and Shen, Jie and Lin, Guang},
  journal={Journal of Computational Physics},
  volume={498},
  pages={112638},
  year={2024},
  publisher={Elsevier}
}

@article{schwerdtner2024nonlinear,
  title={Nonlinear embeddings for conserving {Hamiltonians} and other quantities with Neural {Galerkin} schemes},
  author={Schwerdtner, Paul and Schulze, Philipp and Berman, Jules and Peherstorfer, Benjamin},
  journal={SIAM Journal on Scientific Computing},
  volume={46},
  number={5},
  pages={C583--C607},
  year={2024},
  publisher={SIAM}
}

@article{tang2026energy,
  title={Energy Dissipation Preserving Feature-based {DNN} {Galerkin} Methods for Gradient Flows},
  author={Tang, Tao and Yang, Jiang and Zhao, Yuxiang and Zhu, Quanhui},
  journal={arXiv preprint arXiv:2603.14029},
  year={2026}
}

@article{shi2026preserving,
  title={Preserving conservation laws in the time-evolving natural {Galerkin} method via relaxation and projection},
  author={Shi, Zihao and Wang, Dongling},
  journal={Communications in Nonlinear Science and Numerical Simulation},
  pages={110743},
  year={2026},
  publisher={Elsevier}
}

@article{hu2024energetic,
  title={Energetic variational neural network discretizations of gradient flows},
  author={Hu, Ziqing and Liu, Chun and Wang, Yiwei and Xu, Zhiliang},
  journal={SIAM Journal on Scientific Computing},
  volume={46},
  number={4},
  pages={A2528--A2556},
  year={2024},
  publisher={SIAM}
}

@article{qiao2026unified,
  title={A unified variational framework and lattice {Boltzmann} method for non-isothermal two-phase flows with phase transitions and wetting effects},
  author={Qiao, Zhonghua and Xu, Xianmin and Yang, Xuguang and Zhang, Yuze},
  journal={Journal of Fluid Mechanics},
  volume={1029},
  pages={A30},
  year={2026},
  publisher={Cambridge University Press}
}

@article{yu2021onsagernet,
  title={{OnsagerNet}: Learning stable and interpretable dynamics using a generalized {Onsager} principle},
  author={Yu, Haijun and Tian, Xinyuan and E, Weinan and Li, Qianxiao},
  journal={Physical Review Fluids},
  volume={6},
  number={11},
  pages={114402},
  year={2021},
  publisher={APS}
}

@article{huang2022variational,
  title={Variational {Onsager} Neural Networks {(VONNs)}: A thermodynamics-based variational learning strategy for non-equilibrium {PDEs}},
  author={Huang, Shenglin and He, Zequn and Reina, Celia},
  journal={Journal of the Mechanics and Physics of Solids},
  volume={163},
  pages={104856},
  year={2022},
  publisher={Elsevier}
}

@article{chang2026unsupervised,
  title={Unsupervised operator learning approach for dissipative equations via {Onsager} principle},
  author={Chang, Zhipeng and Wen, Zhenye and Zhao, Xiaofei},
  journal={SIAM Journal on Scientific Computing},
  volume={48},
  number={4},
  pages={C1060--C1085},
  year={2026},
  publisher={SIAM}
}

@article{doi2011onsager,
  title={Onsager's variational principle in soft matter},
  author={Doi, Masao},
  journal={Journal of Physics: Condensed Matter},
  volume={23},
  number={28},
  pages={284118},
  year={2011}
}
\appendix
\section{Supplementary Numerical Details}
\label{app:experimental_details}

This appendix collects the implementation settings, sampling protocols,
reference-solution construction, and supplementary diagnostics used in
Section~4. The definitions, principal accuracy tables, and main
stability observations are retained in the main text.

\subsection{Implementation Environment and Error Evaluation}
\label{app:common_setup}

All numerical experiments are implemented in Python using JAX and are
performed on a single NVIDIA H200 NVL GPU with approximately $140$ GB of
memory. The computing server is equipped with two Intel Xeon 6530P
processors and runs Ubuntu 24.04.5 LTS. The software environment uses
Python 3.11.16 and JAX 0.10.2 with GPU acceleration through the CUDA 13
backend. Unless otherwise stated, all computations are carried out in
double precision.

For weighted formulations, such as the density-ratio representation
$\rho=u\rho_\infty$ used for the Fokker--Planck equation, $M_h(t)$ denotes
the corresponding physical mass and is evaluated as
\begin{equation}
\label{eq:weighted_numerical_mass}
    M_h(t)
    =
    \int_{\Omega}
    u_h(x,t)\rho_\infty(x)\,dx.
\end{equation}
For the computational-cost comparisons, we report the average wall-clock
time per physical time step. The reported solver times exclude the initial
network fitting and JIT-compilation overhead, so that they primarily reflect
the cost of the time-stepping procedures themselves. All methods within 
each comparison are evaluated on the same GPU under identical hardware conditions.

\subsection{One-Dimensional Benchmark Settings and Optimization Cost}
\label{app:one_d_settings}

\paragraph{One-dimensional heat equation}
For all methods, we use the same spatial discretization and neural-network architecture whenever 
applicable. The computational domain $[-\pi,\pi]$ is represented by $N_x=1024$ 
uniformly distributed points. 
Periodicity is incorporated through the Fourier
feature map $x\mapsto[\sin(x),\cos(x)]$,
which is applied to the state network in EDNN
and to the flux networks in the ONTS methods. 
The neural network consists of four hidden layers, each containing $20$ neurons, with the sine 
activation function, giving the architecture $[2, 20, 20, 20, 20, 1]$. 
For the initialization stage, Adam is used for $20\,000$ 
epochs with a batch size of $128$, an initial learning rate of $10^{-3}$, and a weight-decay 
parameter of $10^{-6}$. The same network architecture and initialization
hyperparameters are used for all four methods,
while the quantities represented by the networks
and the subsequent time-stepping procedures differ.

The computational cost per time step is also reported in
Table~\ref{tab:heat_accuracy}. The dominant cost of EDNN and E-ONTS
comes from solving the least-squares problem arising from the evolution
equation, whereas I-ONTS mainly incurs the cost of neural-network
optimization. In the present experiments, the neural-network optimization
in I-ONTS is performed for 40 iterations at each time step. Since this
one-dimensional heat equation is relatively simple, a small number of
optimization iterations is already sufficient to obtain accurate solutions.
The average per-step costs of EDNN, E-ONTS, and I-ONTS are approximately
$0.094$, $0.121$, and $0.059$ seconds, respectively, and are therefore of
the same order of magnitude. 
In contrast, PC-ONTS combines an explicit predictor based on a
least-squares solve with an implicit neural-network correction, where
40 optimization iterations are also used at each time step. It therefore
incurs both types of computational cost, resulting in an average cost of
approximately $0.154$ seconds per time step. Overall, the computational
costs of the four methods remain of the same order for this simple problem,
while the somewhat larger cost of PC-ONTS mainly reflects its combined
predictor--corrector structure. In subsequent experiments, we further
investigate the effect of varying the number of neural-network optimization
iterations at each time step, in order to assess whether the predictor stage
in PC-ONTS can accelerate the convergence of the implicit correction and
thereby reduce the overall computational cost.

\paragraph{One-dimensional Allen--Cahn benchmarks}
The spatial domain is represented by $N_x=1024$ uniformly
distributed points. Periodicity is
incorporated through Fourier features adapted to the corresponding spatial
domain. The neural networks contain four hidden layers with $20$ neurons per
layer and use the sine activation function. 
For the one-dimensional Allen--Cahn benchmarks, EDNN directly 
represents the state $u_\theta$, whereas the Onsager-based methods 
represent the temporal rate. 
The initial neural representation is fitted using the same training 
settings as in Section~\ref{experiment_1dheat}.
Since closed-form transient solutions are unavailable for the following
examples, precomputed high-accuracy numerical solutions are used as references.

\subsection{Two-Dimensional Stiff-Benchmark Implementation}
\label{app:two_d_settings}

\paragraph{Heat equation with varying diffusivity}
The network architectures, collocation points, and optimization 
hyperparameters are kept fixed for all tested values of $\kappa$.
The initial state is independent of $\kappa$, 
whereas the initial target flux $\bm{j}^0=-\kappa \nabla u_0$ 
is fitted separately for each diffusion coefficient. 
We use a $32\times32$ uniform periodic collocation grid, with
the Fourier feature embedding
$(x,y)\mapsto[\sin x,\cos x,\sin y,\cos y]$.
EDNN uses a scalar state network with four hidden layers of width $30$,
whereas the ONTS methods use two independent flux-component networks,
each with four hidden layers of width $20$. All networks employ the sine
activation function. The initial state (for EDNN) or initial flux
(for ONTS) is fitted using Adam for $20{,}000$ epochs with a batch
size of $128$, an initial learning rate of $10^{-3}$, and a weight-decay
coefficient of $10^{-6}$. For I-ONTS and PC-ONTS, the implicit
minimization uses $200$ Adam iterations per physical time step with a
learning rate of $10^{-2}$. The initial state is independent of
$\kappa$, whereas the initial target flux
$\bm j^0=-\kappa\nabla u_0$ is fitted separately for each $\kappa$.

\paragraph{Two-dimensional Cahn--Hilliard problem}
For the numerical comparison, all methods use a
$64\times64$ uniform periodic collocation grid. 
Periodicity is imposed through the
Fourier feature embedding
$(x,y)\mapsto[\sin(\pi x),\cos(\pi x),\sin(\pi y),\cos(\pi y)]$.
EDNN uses one scalar state network with four hidden layers of width $30$,
while the ONTS methods use two independent flux-component networks,
each with four hidden layers of width $20$. All networks employ the sine
activation function. Initial network fitting uses Adam with
$20{,}000$ epochs, a batch size of $256$, and a learning rate of
$10^{-3}$. The implicit I-ONTS and PC-ONTS updates use up to $200$
Adam iterations per time step with a learning rate of $2\times10^{-4}$.

Compared with the second-order heat equation, the Cahn--Hilliard
equation introduces a much stronger high-frequency stiffness through the
fourth-order term $-\lambda\epsilon^2\Delta^2\phi$.
For a Fourier mode with wave number $\mathbf k$, this contribution scales
as $-\lambda\epsilon^2|\mathbf k|^4$, making an explicit temporal
treatment increasingly restrictive for high-frequency components.
The initial two-bubble profile is formed by taking the pointwise
maximum of two diffuse-interface functions and therefore has reduced
regularity near their junction. Since the Cahn--Hilliard flux
$\bm j=-\lambda\nabla\mu$ involves third-order spatial derivatives,
directly fitting the initial flux $\bm j^0$ can be sensitive to the
resulting high-frequency components. To obtain a more regular initial
representation, we instead fit the neural flux networks to 
the reference flux $\bm j^1$ at the first time level $t=\Delta t=10^{-3}$, 
extracted from
the Fourier pseudospectral reference trajectory. This
reference-informed initialization is used only at the initial stage;
subsequent time steps are advanced by the respective ONTS schemes
without using later reference fluxes.

A reference solution is computed using a Fourier pseudo-spectral
method on a $256\times256$ grid with a time step $\Delta t_{\rm ref}=10^{-4}$,
treating the fourth-order linear term implicitly and the
nonlinear term explicitly. 
The reference fields are then spectrally downsampled to a
$128\times128$ reference dataset and further to the $64\times64$
collocation grid used for the neural-solver comparisons. 
Figure~\ref{fig:ch2d_reference_phase} presents snapshots
of the reference solution, illustrating the evolution of the phase-field morphology.
\begin{figure}[H]
    \centering

    \begin{subfigure}{0.98\textwidth}
        \centering
        \includegraphics[width=\textwidth]
        {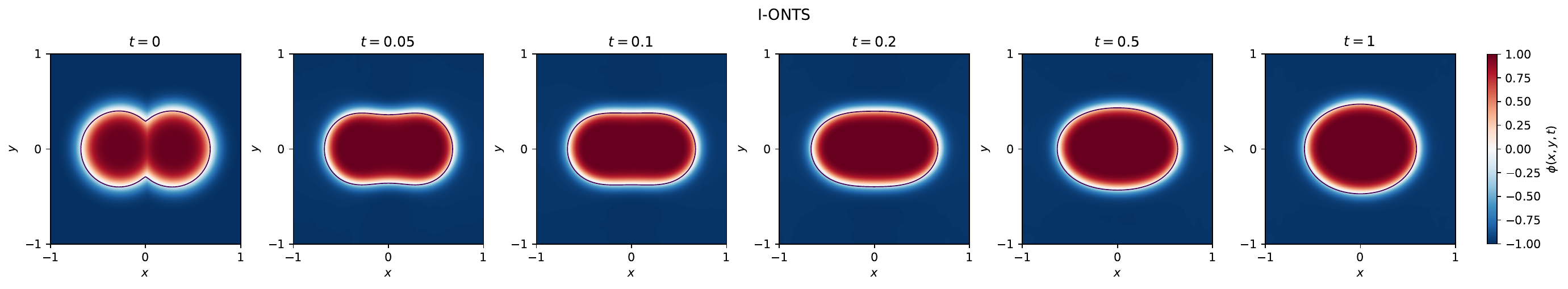}
        \caption{Phase-field evolution obtained by I-ONTS.}
        \label{fig:ch2d_ionts_phase}
    \end{subfigure}

    \vspace{0.4em}

    \begin{subfigure}{0.98\textwidth}
        \centering
        \includegraphics[width=\textwidth]
        {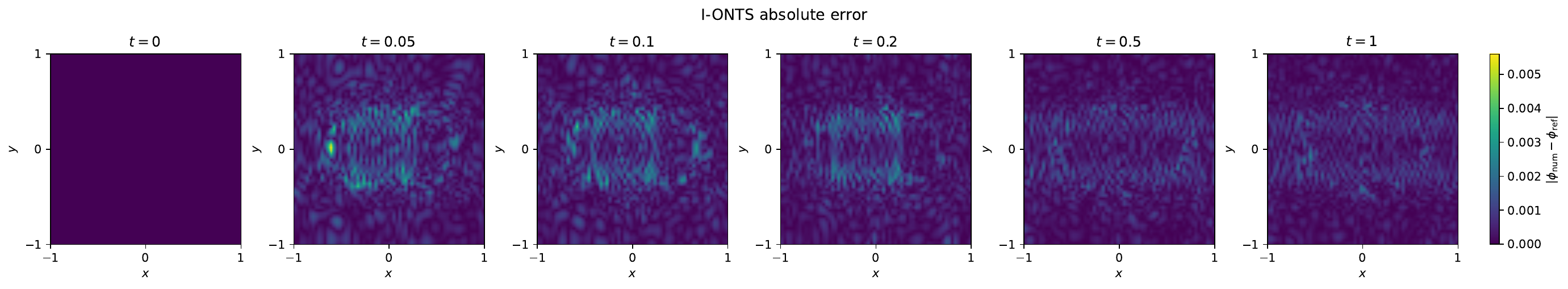}
        \caption{
        Pointwise absolute error
        $|\phi_{\mathrm{I\text{-}ONTS}}-\phi_{\mathrm{ref}}|$.
        }
        \label{fig:ch2d_ionts_error}
    \end{subfigure}

    \caption{
    Phase-field evolution and pointwise absolute error of I-ONTS
    for the two-dimensional Cahn--Hilliard problem.
    I-ONTS remains stable throughout the evolution and
    accurately captures the relaxation of the initially
    connected two-bubble configuration.
    }
    \label{fig:ch2d_ionts_phase_error}
\end{figure}

\begin{figure}[H]
    \centering

    \begin{subfigure}{0.98\textwidth}
        \centering
        \includegraphics[width=\textwidth]
        {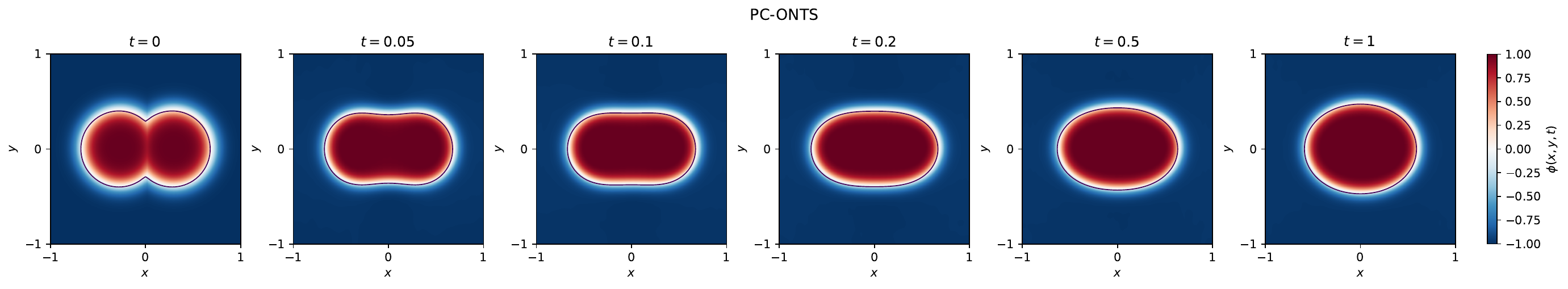}
        \caption{Phase-field evolution obtained by PC-ONTS.}
        \label{fig:ch2d_pconts_phase}
    \end{subfigure}

    \vspace{0.4em}

    \begin{subfigure}{0.98\textwidth}
        \centering
        \includegraphics[width=\textwidth]
        {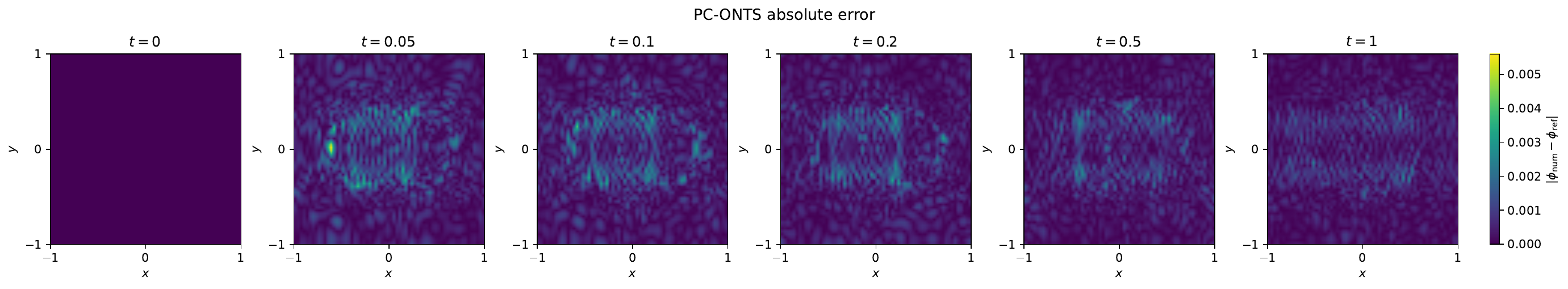}
        \caption{
        Pointwise absolute error
        $|\phi_{\mathrm{PC\text{-}ONTS}}-\phi_{\mathrm{ref}}|$.
        }
        \label{fig:ch2d_pconts_error}
    \end{subfigure}

    \caption{
    Phase-field evolution and pointwise absolute error of PC-ONTS
    for the two-dimensional Cahn--Hilliard problem.
    PC-ONTS also remains stable and closely reproduces the reference
    interfacial dynamics.
    }
    \label{fig:ch2d_pconts_phase_error}
\end{figure}

\paragraph{Supplementary Cahn--Hilliard Phase-Field Comparisons}
Figure \ref{fig:ch2d_ionts_phase_error} and Figure \ref{fig:ch2d_pconts_phase_error} provide additional comparisons of the
phase-field evolution and pointwise absolute errors for the two
implicit Onsager schemes.

\subsection{High-Dimensional Benchmark Implementation}
\label{app:highdim_settings}

\paragraph{Twenty-dimensional Allen--Cahn experiments}
All methods use the same periodic neural representation based on the
features
\begin{equation*}
    \bigl[
        \sin(x_1),\ldots,\sin(x_d),
        \cos(x_1),\ldots,\cos(x_d)
    \bigr],
\end{equation*}
so that periodicity is imposed exactly by construction. 
All neural networks consist of three hidden layers of width $32$ with
the $\tanh$ activation function, and all computations are performed in
double precision. Before the time evolution, the initial neural representations are fitted
using $20{,}000$ training samples and $5{,}000$ validation samples. The
initial fitting is performed for $3000$ Adam epochs with a batch size of
$512$, a learning rate of $10^{-3}$, and a weight-decay coefficient of
$10^{-10}$. The same training configuration is used to initialize the rate networks in
all Onsager-based methods. At each stage of the time evolution, 
spatial quantities are evaluated using
$4096$ Monte Carlo collocation points sampled uniformly from
$[-\pi,\pi]^d$. The same set of collocation points is used for all methods
to ensure a consistent comparison. The reported relative $L^2$ errors are
evaluated independently using an additional set of $20{,}000$ Monte Carlo
test points, also sampled uniformly from $[-\pi,\pi]^d$, which are not used
in the time evolution or optimization. For I-ONTS and PC-ONTS, the implicit
Onsager minimization at each physical time step is carried out using Adam
with learning rate $10^{-4}$ and $200$ neural-network optimization
iterations.

\paragraph{Fokker--Planck Reference Construction}
For the anisotropic Fokker--Planck example, the auxiliary function
$h$ appearing in \eqref{eq:highdim_fp_reference} solves
\begin{equation}
\label{eq:highdim_fp_reference_1d}
    \partial_s h
    =\partial_{xx}h-V'(x)\partial_x h,
    \qquad
    h(x,0)=\cos(2\pi x)-c_\alpha.
\end{equation}
Since the differential operator is separable and the initial
perturbation is additive, the representation
\eqref{eq:highdim_fp_reference} follows from superposition.
The one-dimensional problem is resolved by a Fourier spectral method
with $256$ modes, avoiding a high-dimensional reference grid.

\paragraph{Fokker--Planck sampling and networks}
The spatial integrals are evaluated by Monte Carlo sampling
from the analytically known equilibrium density $\rho_\infty$.
Specifically,
\begin{equation*}
    \int_{\mathbb T^d}
    f(\bm{x})\rho_\infty(\bm{x})\,\mathrm{d}\bm{x}
    \approx
    \frac{1}{N}
    \sum_{j=1}^{N} f(\bm{x}^{(j)}),
    \qquad
    \bm{x}^{(j)}\sim\rho_\infty.
\end{equation*}
This sampling strategy avoids tensor-product grids
and concentrates samples in regions of high probability mass. 
All four methods use the same set of $4096$ samples from 
$\rho_\infty$ for the common spatial collocation and 
least-squares evaluations. For I-ONTS and PC-ONTS, 
the implicit optimization additionally uses 
stochastic mini-batches drawn from a persistent pool
of $32{,}768$ samples, with the corresponding state values
updated during the evolution.
At each Adam iteration, a fresh mini-batch of $4096$
points is sampled with replacement from this pool. 
A separate set of $8192$ samples is reserved for
merit evaluation and candidate selection, independently
of the current optimization mini-batch.

Each component of the rescaled flux is represented
by an independent neural network with four hidden
layers of width $8$, while EDNN employs a scalar
state network with four hidden layers of width $32$.
All networks use sinusoidal activation functions.
To evaluate numerical accuracy independently of the
training and optimization samples, we use an additional
fixed test set of $20{,}000$ points sampled from
$\rho_\infty$.
These points are excluded from time evolution,
least-squares updates, implicit optimization, and
merit-based candidate selection.
The same test set is used for all four methods,
and the resulting Monte Carlo estimates of the
space-time relative $L^2(\rho_\infty)$ errors are
reported in Table~\ref{tab:highdim_fp_accuracy}.

\subsection{Multi-Species Poisson--Nernst--Planck Implementation}
\label{app:pnp_settings}

All methods use a fixed $64\times64$ periodic collocation grid. 
The Poisson equation is solved using a Fourier spectral method, 
with the zero Fourier mode of the charge density removed to satisfy 
the periodic solvability condition. The zero Fourier mode of the 
electrostatic potential is set to zero to fix the gauge. 
To enforce periodicity in the neural representations, the inputs are
mapped to the features $(\sin x,\sin y,\cos x,\cos y)$.
EDNN employs one scalar concentration network per species, each with
four hidden layers of width $32$, whereas the ONTS methods employ six
independent scalar flux networks, each with three hidden layers of
width $16$. All networks use the $\tanh$ activation function.
The EDNN concentration networks and the initial ONTS flux networks
are fitted using $12{,}000$ uniformly sampled training points,
with an additional $3{,}000$ points used for validation.
For ONTS, the initial concentrations are evaluated directly from
\eqref{eq:pnp3_initial}. In the implicit schemes,
$200$ Adam iterations with learning rate $2\times10^{-4}$ are used
per physical time step. The numerical objective also includes a 
small positivity floor and a quadratic penalty for concentrations 
falling below this threshold.

Since no closed-form transient solution is available, we construct a
numerical reference using a Fourier pseudospectral discretization
on a separate $96\times96$ grid, together with the classical
fourth-order Runge--Kutta method with a time step of $10^{-4}$. 
The numerical solutions are evaluated on the same finer grid, 
independently of the $64\times64$ collocation grid used for 
neural time evolution, without feeding the resulting errors 
back into the neural updates. 
The space-time relative $L^2$ error is defined by
\begin{equation}
\label{eq:pnp3_spacetime_error}
    \mathrm{Err}_{\mathrm{st}}
    =\left(
        \frac{
            \displaystyle\sum_{k=0}^{N_t}
            \sum_{i=1}^3
            \bigl\|c_{i,\mathrm{num}}(\cdot,t_k)
                  -c_{i,\mathrm{ref}}(\cdot,t_k)\bigr\|_{L^2(\Omega)}^2
        }{
            \displaystyle\sum_{k=0}^{N_t}
            \sum_{i=1}^3
            \bigl\|c_{i,\mathrm{ref}}(\cdot,t_k)\bigr\|_{L^2(\Omega)}^2
        }
    \right)^{1/2},
\end{equation}
where the spatial norms are approximated using the 
trapezoidal rule on the finer evaluation grid.

\end{document}